\documentclass[11pt]{amsart}

\usepackage[T1]{fontenc}

\usepackage{geometry}
\usepackage{amsmath}
\usepackage{amssymb}
\usepackage[foot]{amsaddr}
\usepackage[dvipsnames]{xcolor}
\usepackage{enumitem}
\usepackage{amscd}
\usepackage{mathrsfs}
\usepackage{cases}
\usepackage[active]{srcltx}
\usepackage[colorlinks,linkcolor={blue},
citecolor={blue},urlcolor={red},]{hyperref}
\usepackage{hyperref}

\newcommand{\R}{\mathbb{R}}

\newcommand{\PP}{\mathbb{P}}

\newcommand{\E}{\mathbb{E}}

\def\qed{\hfill \hbox{$\Box$} \smallskip}
\def\bal{\[\begin{aligned}}
\def\eal{\end{aligned}\]}

\calclayout

\newtheorem{theorem}{Theorem}[section]
\newtheorem{lemma}[theorem]{Lemma}

\newtheorem{assumption}[theorem]{Assumption}
\newtheorem{corollary}[theorem]{Corollary}

\newtheorem{remark}[theorem]{Remark}

\numberwithin{equation}{section}

\begin{document}  

\title[Multiscale Optimal Control with Multiplicative Noise]{Multiscale Linear-Quadratic Stochastic Optimal Control with Multiplicative Noise}

\author{Beniamin Goldys$^*$}
\author{Gianmario Tessitore$^\dagger$}
\author{James Yang$^*$}
\author{Zhou Zhou$^*$}

\address{$^*$ School of Mathematics and Statistics, University of Sydney, Australia. }
\address{$^\dagger$Dipartimento di Matematica e Applicazioni, Universita Milano-Bicocca, Italia.}

\date{\today}
\maketitle

\begin{abstract}

We investigate the asymptotic properties of a finite-time horizon linear-quadratic optimal control problem driven by a multiscale stochastic process with multiplicative Brownian noise. We approach the problem by considering the associated differential Riccati equation and reformulating it as a classical and deterministic singular perturbation problem. Asymptotic properties of this deterministic problem can be gathered from the well-known Tikhonov Theorem. Consequently, we are able to propose two approximation methods to the value function of the stochastic optimal control problem. The first is by constructing an approximately optimal control process whilst the second is by finding the direct limit to the value function. Both approximation methods rely on the existence of a solution to a coupled differential-algebraic Riccati equation with certain stability properties \--- this is the main difficulty of the paper.

\end{abstract}

\section{Introduction}

In this paper, we study a class of linear-quadratic optimal control problems on a finite time interval $[0,T]$. The main assumption is that the dynamical system is comprised of a slow process $X_1$ and a fast process $X^\epsilon_2$ satisfying the following linear stochastic differential equations with multiplicative noise
\begin{align*}
\begin{cases}
dX_1(t) = \left[ A_{11} X_1(t) + A_{12} X^\epsilon_2(t) + B_1 u(t) \right] dt + \left[ C_{11} X_1(t)+ C_{12} X^\epsilon_2(t) + D_1 u(t)\right] dW(t),\\
dX^\epsilon_2(t) = \frac{1}{\epsilon} \left[ A_{21} X_1(t) + A_{22} X^\epsilon_2(t) + B_2 u(t) \right] dt + \frac{1}{\sqrt{\epsilon}} \left[ C_{21} X_1(t)+ C_{22} X^\epsilon_2(t) + D_2 u(t)\right] dW(t),\\
X_1(0) = x_1,\ X^\epsilon_2(0) = x_2,
\end{cases}
\end{align*}
where $\epsilon$ is a small positive parameter representing the ratio between the evolutionary speeds of the slow and fast processes. The objective of the optimal control problem is to minimise a quadratic cost functional with respect to the control process $u$. Our interests lie in deriving estimates for the value function to the optimal control problem when $\epsilon$ is small. Naturally, a \textit{reduced} version of the optimal control problem when $\epsilon$ is formally set to be zero should be considered. However, care needs to be exercised in the fast component due to the degenerate nature of the differential term. For this reason, convergence problems of this type are non-trivial and are commonly referred to as \textit{singular perturbation} problems. Suitably, we shall refer to the optimal control problem described in this paper as a \textit{singularly perturbed linear-quadratic stochastic optimal control problem}, or Problem (SLQP) for short.

Singularly perturbed stochastic optimal control problems have been studied under various formulations and assumptions, see for example \cite{alvarez2002viscosity, borkar2007averaging, dragan2005linear, dragan2012linear, guatteri2018singular, kabanov1991optimal, kabanov1996control, kokotovic1968singular, kokotovic1976singular,kushner2012weak,sannuti1969near, singh1982linear, swikechsingular}. Specifically, in the non-linear case, Alvarez and Bardi \cite{alvarez2002viscosity} formulate the value function as the viscosity solution of a Hamilton-Jacobi-Bellman (HJB) equation, whilst assuming the drift, diffusion and cost functions to be periodic in the fast variable. This assumption amongst others is necessary to ensure the stability of the fast component of the system and the solvability of a reduced HJB equation characterised by a so-called effective Hamiltonian. This HJB approach was extended to the infinite dimensional setting by \'{S}wi\k{e}ch \cite{swikechsingular}. An alternate treatment of the non-linear case with infinite dimensional spaces was studied by Guatteri and Tessitore \cite{guatteri2018singular}. In this case, the authors approach the problem by expressing the value function as the solution of a backwards stochastic differential equation (BSDE) under dissipativity stability assumptions in the drift of the fast variable. In the linear-quadratic setting, a common approach is to take a first-order representation of the solution to the associated Riccati equation. This was first done in the deterministic setting by Kokotovic and Sannuti \cite{kokotovic1968singular, sannuti1969near} and more recently, in the infinite time horizon version of Problem (SLQP) by Dragon et al. \cite{dragan2012linear}.

Traditionally, the benefits of studying singular perturbation problems are twofold. First, asymptotic estimates can be made using the solution of the reduced problem when $\epsilon$ is formally set to be zero, which in some cases are more straightforward to solve such as when the state equations are deterministic or have additive noise, see \cite{kokotovic1968singular, sannuti1969near}. Second, the reduced problem is of a lower order of dimensionality and thus reduces the complexity of the original singular perturbation problem. Applications of the deterministic singularly perturbed optimal control problems has seen active research across a variety of disciplines such as aerospace engineering, biology and chemistry, see the extensive surveys by Naidu \cite{naidu1988singular, naidu2002singular} and references therein. In the stochastic setting, singularly perturbed optimal control have been recently applied to problems in filtering for optimal control problems in finance with partial information by Fouque et al. \cite{fouque2015filtering, fouque2017perturbation} and Kushner \cite{kushner2012weak}. 

It is well-known that for $\epsilon$ fixed, the optimal control and value function of Problem (SLQP) can be characterised in terms of the solution to the associated Riccati equation, see \cite{yong1999stochastic}. By writing the unique solution of the Riccati equation in a first-order representation \cite{dragan2012linear,kokotovic1999singular, xu1997infinite}, the Riccati equation can be rewritten in terms of a system of ODEs, which we call the \textit{full system}. We observe that the full system is in the form of a classical and deterministic singular perturbation problem. Moreover, when $\epsilon$ is formally set to $0$, we obtain a so-called \textit{reduced system} that can be shown to be equivalent to a coupled differential-algebraic Riccati equation of reduced dimensionality. The asymptotic relationship between solutions to the full and reduced systems is non-trivial as the solution of the reduced system is unable to satisfy all the boundary conditions prescribed by the full system. Using the theory of Tikhonov \cite{khalil1996nonlinear, tikhonov1952systems}, the gap in the asymptotic relationship between the full and reduced systems is shown to be resolved by the boundary-layer problem. This is the centrepiece of the paper. From this result, we propose two approaches to estimate the value function of the Problem (SLQP). The first one is by constructing an approximately optimal control process based on the solution to the reduced system. The second one is by directly applying the Tikhonov results to obtain the limiting value function. Both approaches give estimates of the value function with an error of order $O(\epsilon)$.

To the best of our knowledge, Problem (SLQP) on a finite time horizon has not been studied in literature. In the non-linear case \cite{alvarez2002viscosity, guatteri2018singular}, the cost function is assumed to be uniformly Lipschitz in the state variables and thus excludes the quadratic case, as in our paper. Our approach adopts the same first-order representation outlined in the infinite time version \cite{dragan2012linear} and deterministic case \cite{kokotovic1968singular, sannuti1969near}. For this approach to work, we require the existence of the solution to the reduced system with stability properties. It can be shown that the reduced system is equivalent to a pair of so-called reduced differential-algebraic Riccati equations (DARE). In the case of the deterministic case \cite{kokotovic1968singular, sannuti1969near}, which can be easily extended to include additive noise, the reduced DARE is decoupled and thus, the solvability is well-known. However in our case, the multiplicative noise leads to a coupled reduced DARE. The solvability of such a coupled DARE is the main difficulty of this paper and, to the best of our understanding, has not been studied in literature. For this reason, we regard this as a significant contribution of this paper.

The paper proceeds as follows: In Section \ref{section: formulation}, we formulate rigorously Problem (SLQP), the associated Riccati equation and the optimality results. In Section \ref{section: systens}, we partition the solution of the Riccati equation using the first-order representation and derive the full and reduced systems. Sections \ref{section: equivalence} and \ref{section: reducedDARE} are dedicated to showing that the reduced system is equivalent to a coupled reduced DARE and the existence of a stabilising solution. In Section \ref{section: tikhonov}, we tie together the previous three sections via the Tikhonov theorem and establish the convergence properties of the solution to the Riccati equation. Finally, in Section \ref{section: estimates} we formulate an approximating optimal feedback control process based on the solution of the reduced system as well as the limiting value function.

\section{Mathematical Formulation}\label{section: formulation}
\subsection{Notation}

Given a real and separable Hilbert space $E$, the inner product of its elements is denoted by $\langle \cdot, \cdot \rangle_E$ and the associated norm as $|\cdot |_E$. If $G = E\times F$ is the Cartesian product of the Hilbert spaces $E$ and $F$ then $G$ endowed with the inner product $\langle \cdot,\cdot \rangle_G = \langle \cdot,\cdot \rangle_E + \langle \cdot,\cdot \rangle_F$ is also a Hilbert space. For the most part, when there is no confusion, we will drop the subscript in the inner product and norm. For matrices in the space $\R^{n\times n}$, we will write $\mathbb{S}^n$ as the space of symmetric matrices, $\mathbb{S}^n_+$ as the space of non-negative (or positive semi-definite) symmetric matrices and $\mathbb{S}^n_{++}$ as the space of positive (definite) symmetric matrices. For a matrix $M$, we shall use the notation $M^*$ to denote the transpose of $M$ and $M^{-*}$ to denote the inverse of the transpose.

Let $T > 0$ be a finite time horizon and let $W = (W(t))_{0\leq t\leq T}$ be a $1$-dimensional Brownian motion defined on a complete filtered probability space $(\Omega, \mathcal{F}, (\mathcal{F}_t),\PP)$, where $(\mathcal{F}_t)$ is the natural filtration of $W$ augmented by all the $\PP$-null sets in $\mathcal{F}$. We define the following spaces of processes with respect to the Hilbert space $E$:
\begin{itemize}
\item $C([0,T]; E)$ is the space of continuous mappings $F: [0,T] \rightarrow E$ equipped with the norm
\[
\| F\|_{C([0,T]; E)} := \sup_{t\in [0,T]} |F(t)|_{E}.
\] 
\item $L^2_\mathcal{F}(\Omega \times [0,T]; E)$ is the space of equivalence classes of processes $F \in L^2(\Omega \times [0,T]; E)$ admitting a predictable version and equipped with the norm
\[
\| F\|_{L^2_\mathcal{F}(\Omega \times [0,T]; E)} := \left( \E\int_0^T |F(t)|_{E}^2 dt \right)^{1/2}.
\]
For short, we will write $\mathbb{A}_T^n := L^2_\mathcal{F}(\Omega \times [0,T]; \R^n)$.
\item $L^2_\mathcal{F}(\Omega; C([0,T];E))$ is the space of predictable processes $F: \Omega \times [0,T] \rightarrow E$ with continuous paths in $E$ equipped with the norm
\[
\| F\|_{L^2_\mathcal{F}(\Omega; C([0,T];E))} := \left( \E \sup_{t\in [0,T]} |F(t)|_{E}^2  \right)^{1/2}.
\]
For short, we will write $\mathbb{H}_T^n := L^2_\mathcal{F}(\Omega; C([0,T];\R^n))$.
\end{itemize}
\subsection{Problem Formulation}

For arbitrary $x_1 \in \R^{n_1}, x_2\in \R^{n_2}$ and fixed $0 < \epsilon \leq 1$, the processes $X_1 = (X_1(t))_{0\leq t\leq T}$ and $X^\epsilon_2 = (X^\epsilon_2(t))_{0\leq t\leq T}$ are governed by the following linear stochastic differential equations with multiplicative noise
\begin{equation}\label{states}
\begin{split}
\begin{cases}
dX_1(t) = \left[ A_{11} X_1(t) + A_{12} X^\epsilon_2(t) + B_1 u(t) \right] dt \\
\qquad \qquad  + \left[ C_{11} X_1(t)+ C_{12} X^\epsilon_2(t) + D_1 u(t)\right] dW(t),\\
dX^\epsilon_2(t) = \frac{1}{\epsilon} \left[ A_{21} X_1(t) + A_{22} X^\epsilon_2(t) + B_2 u(t) \right] dt  \\
\qquad \qquad  + \frac{1}{\sqrt{\epsilon}} \left[ C_{21} X_1(t)+ C_{22} X^\epsilon_2(t) + D_2 u(t)\right] dW(t),\\
X_1(0) = x_1,\ X^\epsilon_2(0) = x_2,
\end{cases}
\end{split}
\end{equation}
where $u = (u(t))_{0 \leq t\leq T}$ is the control process taking values in $\R^k$ and $A_{ij}, B_i, C_{ij}, D_i$ are deterministic, time-independent matrices of appropriate dimensions. Typically, $X_1$ is referred to as the slow process and $X_2^\epsilon$ is referred to as the fast process due to the presence of $\epsilon$. Let $n = n_1 + n_2$. The above system of equations can be written in under a compact formulation 
\begin{equation}\label{compactstates}
\begin{split}
\begin{cases}
dX^\epsilon(t) = \left[ A^\epsilon X^\epsilon(t) + B^\epsilon u(t) \right] dt + \left[ C^\epsilon X^\epsilon(t) + D^\epsilon u(t)\right] dW(t),\\
X^\epsilon(0) = x := \begin{pmatrix}
x_1\\ x_2
\end{pmatrix},
\end{cases}
\end{split}
\end{equation}
where for all $t\in [0,T]$,
\begin{align*}
\begin{cases}
X^\epsilon(t) = \begin{pmatrix}
X_1(t)\\ X^\epsilon_2(t)
\end{pmatrix}, \quad 
A^\epsilon = \begin{pmatrix}
A_{11} & A_{12}\\
\frac{1}{\epsilon} A_{21} & \frac{1}{\epsilon} A_{22}
\end{pmatrix},\quad
B^\epsilon = \begin{pmatrix}
B_{1}\\ \frac{1}{\epsilon} B_2
\end{pmatrix},\\
C^\epsilon = \begin{pmatrix}
C_{11} & C_{12}\\
\frac{1}{\sqrt{\epsilon}} C_{21} & \frac{1}{\sqrt{\epsilon}} C_{22}
\end{pmatrix},\quad
D^\epsilon = \begin{pmatrix}
D_1\\ \frac{1}{\sqrt{\epsilon}} D_2
\end{pmatrix}.
\end{cases}
\end{align*}

The goal of the optimal control problem is to minimise the cost functional 
\begin{equation}\label{cost}
\begin{split}
J^\epsilon(x;u) &= \frac{1}{2} \E \int_0^T \left[ \langle Q_{11} X_1(t), X_1(t) \rangle + 2\langle Q_{12}^* X_1(t), X^\epsilon_2(t) \rangle \right.\\
&\left. \quad  + \langle Q_{22} X^\epsilon_2(t), X^\epsilon_2(t) \rangle + \langle R u(t), u(t) \rangle \right] dt\\
&=\frac{1}{2} \E \int_0^T \left[ \langle Q X^\epsilon(t), X^\epsilon(t) \rangle + \langle R u(t), u(t) \rangle \right] dt
\end{split}
\end{equation}
with respect to $u$ from the set of admissible controls $\mathbb{A}_T^k := L^2_\mathcal{F}(\Omega \times [0,T]; \R^k)$. Here $Q_{ij}$ and $R$ are matrices of appropriate dimensions and
\begin{align*}
Q = \begin{pmatrix}
Q_{11} & Q_{12}\\
Q_{12}^* & Q_{22}
\end{pmatrix}.
\end{align*}
If $u$ is admissible (i.e. $u \in \mathbb{A}_T^k$) then equation \eqref{compactstates} admits a unique solution $X^\epsilon \in \mathbb{H}_T^n := L^2_\mathcal{F}(\Omega; C([0,T];\R^n))$ for every $0 < \epsilon\leq 1$. Thus, the coupled system of equations \eqref{states} also admits a unique solution $(X_1, X_2^\epsilon) \in \mathbb{H}_T^{n_1} \times \mathbb{H}_T^{n_2}$ for every $0 < \epsilon \leq 1$. Moreover, the cost functional \eqref{cost} is well-defined for all $u \in \mathbb{A}_T^k$. We shall refer to the optimal control problem described by \eqref{compactstates} and \eqref{cost} as the \textit{singularly perturbed stochastic linear-quadratic optimal control problem}, or Problem (SLQP) for short. Moreover, we say that $V^\epsilon(x) = \inf_{u \in \mathbb{A}_T^k} J^\epsilon(x;u)$ is the value function of Problem (SLQP). Finally, the associated Riccati equation is given by 
\begin{equation}\label{riccati}
\begin{split}
\begin{cases}
&\frac{dP^\epsilon}{dt} + (A^\epsilon)^* P^\epsilon + P^\epsilon A^\epsilon + (C^\epsilon)^* P^\epsilon C^\epsilon + Q\\
&\quad - \left[ (B^\epsilon)^* P^\epsilon  + (D^\epsilon)^* P^\epsilon C^\epsilon \right]^*\left[R + (D^\epsilon)^* P^\epsilon D^\epsilon \right]^{-1} \left[ (B^\epsilon)^* P^\epsilon  + (D^\epsilon)^* P^\epsilon C^\epsilon \right] = 0,\\
&R + (D^\epsilon)^* P^\epsilon(t) D^\epsilon > 0, \quad \forall t\in [0,T],\\
&P^\epsilon(T) = 0.
\end{cases}
\end{split}
\end{equation}

Going forward, the main assumptions of the paper are stated in Assumption \ref{mainassumptions}. The first assumption is sufficient for the existence and uniqueness of a solution to the Riccati equation as well as an optimal control. The second assumption is needed in Section \ref{section: systens} to show that the reduced system is equivalent to a reduced differential-algebraic Riccati equation. Finally, the last assumption is sufficient to show that the reduced differential-algebraic Riccati equation admits a unique stablising solution. This is shown in Section \ref{section: reducedDARE}. 
\begin{assumption}\label{mainassumptions}
Suppose the following holds:
\begin{enumerate}
\item The matrices $Q$ and $R$ are strictly positive definite;
\item $A_{22}$ is invertible;
\item $[A_{22},C_{22}]$ is $L^2$-stable. That is, if for any $x_2 \in \R^n$, the solution to 
\begin{equation}
dX_2(t) = A_{22} X_2(t) dt + C_{22}X_2(t) dW(t),\quad X_2(0) = x_2,
\end{equation}
is a continuous $\mathcal{F}$-adapted process satisfying
\begin{align*}
\E \int_0^\infty |X_2(t)|^2 dt < \infty.
\end{align*}
\end{enumerate}
\end{assumption}

\begin{remark}
From Lemma 2.2 of \cite{sun2018stochastic}, $[A_{22},C_{22}]$ being $L^2$-stable is equivalent to the existence of a $Y \in \mathbb{S}^{n_2}_{++}$ such that
\[
A_{22}^* Y + Y A_{22} + C_{22}^* Y C_{22} < 0.
\]
In other words, $A_{22}$ should be negative enough. 

Moreover, it is possible to weaken this assumption to $[A_{22} + B_2 L, C_{22} + B_2 L]$ being $L^2$-stable for some matrix $L$ of appropriate dimensions. Generally, this is the standard assumption for the fast equation or infinite-time horizon problems, see for example \cite{dragan2012linear, sun2018stochastic}. We avoid this generality in order to reduce some notational complexity but note that the subsequent results will hold albeit with different coefficient matrices. 
\end{remark}

\begin{lemma}\label{lemma: riccati}
Suppose that Assumption \ref{mainassumptions}-(1) holds. Then, for every $0 < \epsilon\leq 1$, the Riccati equation \eqref{riccati} admits a unique solution $P^\epsilon \in C([0,T]; \mathbb{S}^n_+)$.
\end{lemma}
\noindent \textit{Proof.} See Theorem 7.2 in Chapter 6 of \cite{yong1999stochastic}.
\qed

Next, we state the well-known optimality result, which characterises the optimal control and value function in terms of the unique positive semi-definite solution to the Riccati equation.
\begin{theorem}\label{theorem: optimality}
Suppose that Assumption \ref{mainassumptions}-(1) holds. Let $P^\epsilon \in C([0,T];\mathbb{S}^n_+)$ be the unique solution to the Riccati equation \eqref{riccati}. Then, for every $0 < \epsilon \leq 1$, Problem (SLQP) admits an unique optimal control $\widehat{u}^\epsilon \in \mathbb{A}_T^k$ given by
\begin{equation}
\widehat{u}^\epsilon(t) =  \widehat{F}^\epsilon(t) \widehat{X}^\epsilon(t) ,\quad \forall t\in [0,T],
\end{equation}
where $\widehat{X}^\epsilon(t) = X^\epsilon(t; \widehat{u}^\epsilon(t))$ and the feedback operator $\widehat{F}^\epsilon$ is defined by
\begin{equation}\label{fboperator}
\widehat{F}^\epsilon(t) = - (R + (D^\epsilon)^* P^\epsilon(t) D^\epsilon)^{-1} \left[ (B^\epsilon)^* P^\epsilon(t) + (D^\epsilon)^* P^\epsilon(t) C^\epsilon\right], \quad \forall t\in [0,T].
\end{equation}
Moreover, the value function is given by
\begin{equation}
V^\epsilon(x) = \frac{1}{2}\langle P^\epsilon(0)x,x\rangle.
\end{equation}
\end{theorem}
\noindent \textbf{Proof.} See Theorem 6.1 in Chapter 6 of \cite{yong1999stochastic}. \qed

\section{Singular perturbation of the Riccati equation}\label{section: systens}

The goal of this section is to reformulate the Riccati equation \eqref{riccati} as a classical and deterministic singular perturbation problem. Let us partition the solution $P^\epsilon$ to the Riccati equation \eqref{riccati} in the first-order form
\begin{equation}\label{partition}
P^\epsilon(t) = \begin{pmatrix}
P^\epsilon_{11}(t) & \epsilon P^\epsilon_{12}(t)\\
\epsilon (P^\epsilon_{12}(t))^* & \epsilon P^\epsilon_{22}(t)
\end{pmatrix}, \quad \forall t\in [0,T].
\end{equation}
Moreover, let 
\begin{equation}
\begin{split}
\Delta^\epsilon(t) &= R + (D^\epsilon)^* P^\epsilon(t) D^\epsilon \\
&= R + D_1^* P_{11}^\epsilon(t) D_1 + \sqrt{\epsilon} [D_2^* (P^\epsilon_{12}(t))^* D_1 + D_1^* P^\epsilon_{12}(t) D_2] + D_2^* P^\epsilon_{22}(t) D_2.
\end{split}
\end{equation} 
Applying the partition \eqref{partition} to the Riccati equation \eqref{riccati}, we obtain a system of differential equations, which we call the \textit{full system}
\begin{subnumcases}{\label{fullsystem}}
\frac{dP^\epsilon_{11}}{dt} + f(P^\epsilon_{11}, P^\epsilon_{12}, P^\epsilon_{22}, \epsilon)  = 0, \quad  P^\epsilon_{11}(T) = 0,  \\
\epsilon\frac{dP^\epsilon_{12}}{dt} + g_1(P^\epsilon_{11}, P^\epsilon_{12}, P^\epsilon_{22},\epsilon)=0 ,\quad P^\epsilon_{12}(T) = 0,\\
\epsilon\frac{dP^\epsilon_{22}}{dt} + g_2(P^\epsilon_{11}, P^\epsilon_{12}, P^\epsilon_{22}, \epsilon)=0, \quad P^\epsilon_{22}(T) = 0,\\
\Delta^\epsilon(t) > 0 \quad \forall t\in [0,T]\ a.e.,
\end{subnumcases}
where the functions $f, g_1$ and $g_2$ are defined as
\begin{equation}
\begin{split}
&f(P^\epsilon_{11}, P^\epsilon_{12}, P^\epsilon_{22}, \epsilon) \\
&\quad= A^*_{11} P^\epsilon_{11} + P^\epsilon_{11} A_{11} + A_{21}^* (P^\epsilon_{12})^* + P^\epsilon_{12} A_{21} + C_{11}^* P^\epsilon_{11} C_{11} \\
&\qquad + \sqrt{\epsilon} \left( C_{21}^* (P^\epsilon_{12})^* C_{11} + C_{11}^* P^\epsilon_{12} C_{21} \right) + C_{21}^* P^\epsilon_{22} C_{21} + Q_{11} \\
&\qquad - \left[B_1^* P^\epsilon_{11} + B_2^* (P^\epsilon_{12})^* + D_1^* P^\epsilon_{11} C_{11} + \sqrt{\epsilon}\left( D_2^* (P^\epsilon_{12})^* C_{11} + D_1^* P^\epsilon_{12} C_{21} \right) + D_2^* P^\epsilon_{22} C_{21}\right]^* (\Delta^{\epsilon})^{-1} \\
&\qquad \qquad\left[ B_1^* P^\epsilon_{11} + B_2^* (P^\epsilon_{12})^* + D_1^* P^\epsilon_{11} C_{11} + \sqrt{\epsilon}\left( D_2^* (P^\epsilon_{12})^* C_{11} + D_1^* P^\epsilon_{12} C_{21} \right) + D_2^* P^\epsilon_{22} C_{21}\right]
\end{split}
\end{equation}
\begin{equation}
\begin{split}
&g_1(P^\epsilon_{11}, P^\epsilon_{12}, P^\epsilon_{22},\epsilon) \\
&\quad = \epsilon A_{11}^* P^\epsilon_{12} + P^\epsilon_{11} A_{12} + A_{21}^* P^\epsilon_{22} + P^\epsilon_{12} A_{22} + C_{11}^* P^\epsilon_{11} C_{12}\\
&\qquad + \sqrt{\epsilon} \left( C_{21}^* (P^\epsilon_{12})^* C_{12} + C_{11}^* P^\epsilon_{12} C_{22} \right) + C_{21}^* P^\epsilon_{22} C_{22} + Q_{12} \\
&\qquad - \left[ B_1^* P^\epsilon_{11} + B_2^* (P^\epsilon_{12})^* + D_1^* P^\epsilon_{11} C_{11} + \sqrt{\epsilon}\left( D_2^* (P^\epsilon_{12})^* C_{11} + D_1^* P^\epsilon_{12} C_{21} \right) + D_2^* P^\epsilon_{22} C_{21}\right]^* (\Delta^{\epsilon})^{-1} \\
&\qquad \qquad \left[  \epsilon B_1^* P^\epsilon_{12} + B_2^* P^\epsilon_{22} + D_1^* P^\epsilon_{11} C_{12} + \sqrt{\epsilon}\left( D_2^* (P^\epsilon_{12})^* C_{12} + D_1^* P^\epsilon_{12} C_{22} \right) + D_2^* P^\epsilon_{22} C_{22} \right] 
\end{split}
\end{equation}
\begin{equation}
\begin{split}
&g_2(P^\epsilon_{11}, P^\epsilon_{12}, P^\epsilon_{22},\epsilon) \\
&\quad = A_{22}^* P^\epsilon_{22} + P^\epsilon_{22} A_{22} + \epsilon\left(A_{12}^* P^\epsilon_{12} + (P^\epsilon_{12})^* A_{12} \right) + C_{12}^* P^\epsilon_{11} C_{12}\\
&\qquad + \sqrt{\epsilon} \left(C_{22}^* (P^\epsilon_{12})^* C_{12} + C_{12}^* P^\epsilon_{12} C_{22} \right) + C_{22}^* P^\epsilon_{22} C_{22} + Q_{22} \\
&\qquad - \left[ \epsilon B_1^* P^\epsilon_{12} + B_2^* P^\epsilon_{22} + D_1^* P^\epsilon_{11} C_{12} + \sqrt{\epsilon}\left( D_2^* (P^\epsilon_{12})^* C_{12} + D_1^* P^\epsilon_{12} C_{22} \right) + D_2^* P^\epsilon_{22} C_{22}\right]^* (\Delta^{\epsilon})^{-1}\\
&\qquad \qquad \left[  \epsilon B_1^* P^\epsilon_{12} + B_2^* P^\epsilon_{22} + D_1^* P^\epsilon_{11} C_{12} + \sqrt{\epsilon}\left( D_2^* (P^\epsilon_{12})^* C_{12} + D_1^* P^\epsilon_{12} C_{22} \right) + D_2^* P^\epsilon_{22} C_{22}\right].
\end{split}
\end{equation}
It is clear that we have the following lemmas.
\begin{lemma}
For fixed $0 < \epsilon \leq 1$, the full system admits a solution $(P^\epsilon_{11},P^\epsilon_{12},P^\epsilon_{22}) \in C([0,T]; \mathbb{S}^{n_1}_+) 
\times C([0,T];\R^{n_1 \times n_2}) \times C([0,T]; \mathbb{S}^{n_2}_+)$ if and only if the Riccati equation \eqref{riccati} admits a solution $P^\epsilon \in  C([0,T];\mathbb{S}^n_+)$.
\end{lemma}
\begin{lemma}\label{lemma: fullsystemsoln}
Suppose that Assumption \ref{mainassumptions}-(1) holds. Then the full system admits a unique solution $(P^\epsilon_{11},P^\epsilon_{12},P^\epsilon_{22}) \in C([0,T]; \mathbb{S}^{n_1}_+) 
\times C([0,T];\R^{n_1 \times n_2}) \times C([0,T]; \mathbb{S}^{n_2}_+)$.
\end{lemma}

A well-known approach to analysing the asymptotic behaviour of dynamical systems represented in the form of the full system \eqref{fullsystem} as $\epsilon$ tends towards zero is by applying a version of the Tikhonov theorem \cite{khalil1996nonlinear, tikhonov1952systems}. This approach has been used in the context of Riccati equations; both in the differential \cite{kokotovic1968singular,sannuti1969near} and algebraic \cite{dragan2012linear} cases. As an initial step, we need to consider the differential-algebraic system when $\epsilon$ is formally set to be $0$ in the full system. In this case, we obtain the following so-called \textit{reduced system}
\begin{subnumcases}{\label{reducedsystem}}
\frac{d\overline{P}_{11}}{dt} + f(\overline{P}_{11}, \overline{P}_{12}, \overline{P}_{22}, 0) , \quad \overline{P}_{11}(T) = 0,\label{reduced1}\\
g_1(\overline{P}_{11}, \overline{P}_{12}, \overline{P}_{22}, 0)=0,\label{reduced2}\\
g_2(\overline{P}_{11}, \overline{P}_{12}, \overline{P}_{22}, 0) =0,\label{reduced3}\\
\overline{\Delta}(t) > 0, \quad \forall t\in [0,T]\ a.e.,
\end{subnumcases}
where $\overline{\Delta}(t) =R + D_1^* \overline{P}_{11}(t) D_1 + D_2^* \overline{P}_{22}(t) D_2$. 

A prerequisite to applying the Tikhonov theorem is to show that the reduced system admits a solution $(\overline{P}_{11}, \overline{P}_{12}, \overline{P}_{22})$ such that  $(\overline{P}_{12}, \overline{P}_{22})$ is an isolated root of \eqref{reduced2}-\eqref{reduced3} and $\overline{P}_{11}$ is unique relative to this isolated root. We will dedicate the next two sections showing this.

\section{Equivalence to the reduced differential-algebraic Riccati equation}\label{section: equivalence}

In this section, we show that the reduced system \eqref{reducedsystem} is equivalent to a coupled differential-algebraic Riccati equaiton. For $A_{22}$ invertible, let us define the matrices
\begin{equation}\label{reducedmatrices}
\begin{split}
\begin{cases}
A_s = A_{11} - A_{12} A_{22}^{-1} A_{21}, \ B_s = B_1 - A_{12} A_{22}^{-1} B_2,\\
C_{1s} = C_{11} - C_{12} A_{22}^{-1} A_{21}, \ C_{2s} = C_{21} - C_{22} A_{22}^{-1} A_{21},\\
D_{1s} = D_1 - C_{12} A_{22}^{-1} B_2,\ D_{2s} = D_2 - C_{22} A_{22}^{-1} B_2,\\
Q_s = Q_{11} - Q_{12} A_{22}^{-1} A_{21} - A_{21}^* A_{22}^{-*} Q_{12}^* + A_{21}^* A_{22}^{-*} Q_{22} A_{22}^{-1} A_{21},\\
L_s = B_2^* A_{22}^{-*} \left( Q_{22} A_{22}^{-1}A_{21} - Q_{12}^* \right), \ R_s = R + B_2^* A_{22}^{-*} Q_{22} A_{22}^{-1} B_2.
\end{cases}
\end{split}
\end{equation}
Let $\overline{\Delta}_s(t) := R_s + D_{1s}^* \overline{P}_{11}(t) D_{1s} + D_{2s}^* \overline{P}_{22}(t) D_{2s}$. Consider the following coupled differential-algebraic system of equations
\begin{subnumcases}{\label{reducedDARE}}
\label{reducedDRE}
\begin{split}
\frac{d\overline{P}_{11}}{dt} &+ A_s^* \overline{P}_{11} + \overline{P}_{11} A_s + C_{1s}^* \overline{P}_{11} C_{1s} \\
&+ C_{2s}^* \overline{P}_{22} C_{2s} - M_s^* \overline{\Delta}_s^{-1} M_s + Q_s = 0 , \quad \overline{P}_{11}(T) = 0,
\end{split}\\
\label{reducedARE}
A_{22}^* \overline{P}_{22} + \overline{P}_{22} A_{22} + C_{12}^* \overline{P}_{11} C_{12} + C_{22}^* \overline{P}_{22} C_{22} - M_2^* \overline{\Delta}^{-1} M_2 + Q_{22} = 0,\\
\label{convexDRE}
\overline{\Delta}_s(t) > 0,\quad \forall t\in [0,T]\ a.e.,\\
\label{convexARE}
\overline{\Delta}(t) > 0,\quad \forall t\in [0,T]\ a.e.,
\end{subnumcases}
where 
\begin{align*}
\begin{cases}
M_s = B_s^* \overline{P}_{11} + D_{1s}^* \overline{P}_{11} C_{1s} + D_{2s}^* \overline{P}_{22} C_{2s} + L_s,\\
M_2 = B_2^* \overline{P}_{22} + D_1^* \overline{P}_{11} C_{12} + D_2^* \overline{P}_{22} C_{22}.
\end{cases}
\end{align*}
We call the above system the \textit{reduced differential-algebraic Riccati equation} or reduced DARE for short. The existence of a solution to the reduced DARE and consequently the reduced system, is the focus of Section \ref{section: reducedDARE}. 

The equivalence between the reduced system \eqref{reducedsystem} and the reduced DARE \eqref{reducedDARE} follows in the same way as in the infinite time horizon counterpart \cite{dragan2012linear}. We will include this derivation for completeness. It is convenient to begin by introducing the operators $(\overline{F}_1,\overline{F}_2)$ defined as
\begin{subnumcases}{}
\label{reducedoperator1}
\overline{F}_1(t) = - \overline{\Delta}^{-1}(t) \left[ B_1^* \overline{P}_{11}(t) + B_2^* \overline{P}_{12}^*(t) + D_1^* \overline{P}_{11}(t) C_{11} + D_2^* \overline{P}_{22}(t) C_{21} \right],\\
\label{reducedoperator2}
\overline{F}_2(t) = - \overline{\Delta}^{-1}(t) \left[ B_2^* \overline{P}_{22}(t) + D_1^* \overline{P}_{11}(t) C_{12} + D_2^* \overline{P}_{22}(t) C_{22}\right].
\end{subnumcases}
As a result, we can write the reduced system \eqref{reducedsystem} as
\begin{subnumcases}{}
\label{reduced4}
\begin{split}
\frac{d\overline{P}_{11}}{dt} &+ A_{11}^* \overline{P}_{11} + \overline{P}_{11} A_{11} + A_{21}^* \overline{P}_{12}^* + \overline{P}_{12} A_{21} \\
&+ C_{11}^* \overline{P}_{11} C_{11} + C_{21}^* \overline{P}_{22} C_{21} + Q_{11} - \overline{F}_1^* \overline{\Delta} \overline{F}_1 = 0, \quad \overline{P}_{11}(T) = 0,
\end{split}\\
\label{reduced5}
\overline{P}_{11} A_{12} + A_{21}^* \overline{P}_{22} + \overline{P}_{12} A_{22} + C_{11}^* \overline{P}_{11} C_{12} + C_{21}^* \overline{P}_{22} C_{22} + Q_{12} - \overline{F}_1^* \overline{\Delta} \overline{F}_2 = 0,\\
\label{reduced6}
A_{22}^* \overline{P}_{22} + \overline{P}_{22} A_{22} + C_{12}^* \overline{P}_{11} C_{12} + C_{22}^* \overline{P}_{22} C_{22} + Q_{22} - \overline{F}_2^* \overline{\Delta} \overline{F}_2 = 0,\\
\overline{\Delta}(t) > 0,\quad \forall t\in [0,T]\ a.e.
\end{subnumcases}
Before we show the equivalence, we state some invertibility results.
\begin{lemma}\label{lemma: invertible00}
Suppose that Assumption \ref{mainassumptions}-(2) holds and $A_{22} + B_2 \overline{F}_2$ is invertible. Then the matrix $I + \overline{F}_2 A_{22}^{-1} B_2$ is invertible where
\begin{equation}\label{identity1}
(I + \overline{F}_2 A_{22}^{-1} B_2)^{-1} = I - \overline{F}_2 (A_{22} + B_2 \overline{F}_2)^{-1} B_2.
\end{equation}
Furthermore, 
\begin{equation}\label{identity2}
(A_{22} + B_2 \overline{F}_2)^{-1} = A_{22}^{-1} - A_{22}^{-1} B_2 (I + \overline{F}_2 A_{22}^{-1} B_2)^{-1} \overline{F}_2 A_{22}^{-1}. 
\end{equation}
\end{lemma}
\noindent \textit{Proof.} To see that $I - \overline{F}_2 (A_{22} + B_2 \overline{F}_2)^{-1} B_2$ is the inverse of $I + \overline{F}_2 A_{22}^{-1} B_2$, observe that
\begin{align*}
&\left[I + \overline{F}_2 A_{22}^{-1} B_2\right] \left[I - \overline{F}_2 (A_{22} + B_2 \overline{F}_2)^{-1} B_2\right] \\
&\quad = I + \overline{F}_2 \left[ A_{22}^{-1} - (A_{22} + B_2 \overline{F}_2)^{-1} - A_{22}^{-1} B_2 \overline{F}_2 (A_{22} + B_2 \overline{F}_2)^{-1} \right] B_2\\
&\quad = I + \overline{F}_2 \left[ A_{22}^{-1} (A_{22} + B_2 \overline{F}_2) - I - A_{22}^{-1} B_2 \overline{F}_2  \right] (A_{22} + B_2 \overline{F}_2)^{-1} B_2\\
&\quad = I.
\end{align*}
Using \eqref{identity1} we can show \eqref{identity2} as follows
\begin{align*}
A_{22}^{-1} B_2 (I + \overline{F}_2 A_{22}^{-1} B_2)^{-1} \overline{F}_2 A_{22}^{-1} 
&= A_{22}^{-1} B_2 \overline{F}_2 A_{22}^{-1} - A_{22}^{-1} B_2  \overline{F}_2 (A_{22} + B_2 \overline{F}_2)^{-1} B_2 \overline{F}_2 A_{22}^{-1}\\
&= A_{22}^{-1} B_2 \overline{F}_2  A_{22}^{-1} \\
&\quad - A_{22}^{-1} (A_{22} +B_2 \overline{F}_2 - A_{22}) (A_{22} + B_2 \overline{F}_2)^{-1} B_2 \overline{F}_2 A_{22}^{-1}\\
&= (A_{22} + B_2 \overline{F}_2)^{-1} B_2 \overline{F}_2 A_{22}^{-1}\\
&= (A_{22} + B_2 \overline{F}_2)^{-1} (A_{22} + B_2 \overline{F}_2 - A_{22}) A_{22}^{-1}\\
&= A_{22}^{-1} - (A_{22} + B_2 \overline{F}_2)^{-1}.
\end{align*}
\qed

We state the equivalence result below.

\begin{theorem}\label{theorem: equiv}
Suppose that Assumption \ref{mainassumptions}-(2) holds. We have the following:
\begin{enumerate}
\item If $(\overline{P}_{11}, \overline{P}_{12}, \overline{P}_{22})$ is the solution to the reduced system \eqref{reducedsystem} and $A_{22} + B_2 \overline{F}_2$ is invertible, where $\overline{F}_2$ is defined in \eqref{reducedoperator2}, then $(\overline{P}_{11}, \overline{P}_{22})$ is the solution to the reduced differential-algebraic Riccati equation \eqref{reducedDARE}.
\item If $(\overline{P}_{11}, \overline{P}_{22})$ is the solution to the reduced differential-algebraic Riccati equation \eqref{reducedDARE} and $A_{22} + B_2 \overline{F}_2$ is invertible, where $\overline{F}_2$ is defined in \eqref{reducedoperator2}, then $(\overline{P}_{11},\overline{P}_{12},\overline{P}_{22})$ is the solution to the reduced system \eqref{reducedsystem} with
\begin{equation}\label{P12equiv}
\begin{split}
\overline{P}_{12} &=- \left[ \overline{P}_{11} A_{12} + A_{21}^* \overline{P}_{22} + C_{11}^* \overline{P}_{11} C_{12} + C_{21}^* \overline{P}_{22} C_{22} \right.\\
&\left. \qquad + Q_{12} + \left(  \overline{P}_{11}B_1  + C_{11}^* \overline{P}_{11} D_1 + C_{21}^* \overline{P}_{22} D_2 \right) \overline{F}_2 \right] 
\left[ A_{22} + B_2 \overline{F}_2 \right]^{-1}.
\end{split}
\end{equation}
\end{enumerate}
\end{theorem}
\noindent \textit{Proof.} We will start by proving Theorem \ref{theorem: equiv}-(1). Let $(\overline{P}_{11}, \overline{P}_{12}, \overline{P}_{22})$ be the solution to the reduced system \eqref{reducedsystem} and $A_{22} + B_2 \overline{F}_2$ be invertible, where $\overline{F}_2$ is defined in \eqref{reducedoperator2}. Observe that the equation \eqref{reduced3} is precisely the reduced algebraic Riccati equation \eqref{reducedARE} and thus share the same solution $\overline{P}_{22}$. To show that the equation \eqref{reduced1} can be rewritten as the reduced differential Riccati equation \eqref{reducedDRE} we will look to eliminate the $\overline{P}_{12}$ terms in the prior. To this end, we use the invertibility assumption of $A_{22}$ to rewrite \eqref{reduced5} as
\begin{equation}\label{P12variant}
\overline{P}_{12}  = \left[ \overline{F}_1^* \overline{\Delta} \overline{F}_2 - \overline{P}_{11} A_{12} - A_{21}^* \overline{P}_{22} -C_{11}^* \overline{P}_{11} C_{12} - C_{21}^* \overline{P}_{22} C_{22} - Q_{12} \right] A_{22}^{-1}.
\end{equation}
Substituting \eqref{P12variant} into the right-hand side of \eqref{reduced4}, we obtain
\begin{align*}
H &:= A_{11}^* \overline{P}_{11} + \overline{P}_{11} A_{11} + A_{21}^* \overline{P}_{12}^* + \overline{P}_{12} A_{21} \\
&\quad + C_{11}^* \overline{P}_{11} C_{11} + C_{21}^* \overline{P}_{22} C_{21} + Q_{11} - \overline{F}_1^* \overline{\Delta} \overline{F}_1\\
&= A_{11}^* \overline{P}_{11} + \overline{P}_{11} A_{11}  + C_{11}^* \overline{P}_{11} C_{11} + C_{21}^* \overline{P}_{22} C_{21} + Q_{11} - \overline{F}_1^* \overline{\Delta} \overline{F}_1\\
&\quad + A_{21}^* A_{22}^{-*} \left[ \overline{F}_1^* \overline{\Delta} \overline{F}_2 - \overline{P}_{11} A_{12} - A_{21}^* \overline{P}_{22} -C_{11}^* \overline{P}_{11} C_{12} - C_{21}^* \overline{P}_{22} C_{22} - Q_{12} \right]^* \\
&\quad + \left[ \overline{F}_1^* \overline{\Delta} \overline{F}_2 - \overline{P}_{11} A_{12} - A_{21}^* \overline{P}_{22} -C_{11}^* \overline{P}_{11} C_{12} - C_{21}^* \overline{P}_{22} C_{22} - Q_{12} \right] A_{22}^{-1} A_{21}\\
&= A_s^* \overline{P}_{11} + \overline{P}_{11}A_s - A_{21}^* A_{22}^{-*}  \overline{P}_{22} A_{21} - A_{21}^* \overline{P}_{22} A_{22}^{-1} A_{21} \\
&\quad + C_{11}^* \overline{P}_{11} C_{11} - A_{21}^* A_{22}^{-*} C_{12}^* \overline{P}_{11} C_{11} - C_{11}^* \overline{P}_{11} C_{12} A_{22}^{-1} A_{21}\\
&\quad + C_{21}^* \overline{P}_{22} C_{21} - A_{21}^* A_{22}^{-*} C_{22}^* \overline{P}_{22} C_{21} - C_{21}^* \overline{P}_{22} C_{22} A_{22}^{-1} A_{21}\\
&\quad + Q_{11} - A_{21}^* A_{22}^{-*} Q_{12}^* - Q_{12} A_{22}^{-1} A_{21} + A_{21}^* A_{22}^{-*} \overline{F}_2^* \overline{\Delta} \overline{F}_1 + \overline{F}_1^* \overline{\Delta} \overline{F}_2 A_{22}^{-1} A_{21} - \overline{F}_1^* \overline{\Delta} \overline{F}_1.
\end{align*}
Completing the square, we have that
\begin{align*}
H &= A_s^* \overline{P}_{11} + \overline{P}_{11}A_s - A_{21}^* A_{22}^{-*}  \overline{P}_{22} A_{21} - A_{21}^* \overline{P}_{22} A_{22}^{-1} A_{21}\\
&\quad  + C_{1s}^* \overline{P}_{11} C_{1s} + C_{2s}^* \overline{P}_{22} C_{2s} - \left[ \overline{F}_1 - \overline{F}_2 A_{22}^{-1}A_{21}\right]^* \overline{\Delta} \left[ \overline{F}_1 - \overline{F}_2 A_{22}^{-1}A_{21}\right]  \\
&\quad - \left[ C_{12} A_{22}^{-1} A_{21}\right]^* \overline{P}_{11} \left[ C_{12} A_{22}^{-1} A_{21}\right] - \left[ C_{22} A_{22}^{-1} A_{21}\right]^* \overline{P}_{22} \left[ C_{22} A_{22}^{-1} A_{21}\right]\\
&\quad + Q_{11} - A_{21}^* A_{22}^{-*} Q_{12}^* - Q_{12} A_{22}^{-1} A_{21} + A_{21}^* A_{22}^{-*} \overline{F}_2^* \overline{\Delta} \overline{F}_2 A_{22}^{-1} A_{21}.
\end{align*}
From \eqref{reduced6}, we have that
\begin{align*}
H &= A_s^* \overline{P}_{11} + \overline{P}_{11}A_s - A_{21}^* A_{22}^{-*}  \overline{P}_{22} A_{21} - A_{21}^* \overline{P}_{22} A_{22}^{-1} A_{21}\\
&\quad + C_{1s}^* \overline{P}_{11} C_{1s} + C_{2s}^* \overline{P}_{22} C_{2s} - \left[ \overline{F}_1 - \overline{F}_2 A_{22}^{-1}A_{21}\right]^* \overline{\Delta} \left[ \overline{F}_1 - \overline{F}_2 A_{22}^{-1}A_{21}\right]  \\
&\quad - \left[ C_{12} A_{22}^{-1} A_{21}\right]^* \overline{P}_{11} \left[ C_{12} A_{22}^{-1} A_{21}\right] - \left[ C_{22} A_{22}^{-1} A_{21}\right]^* \overline{P}_{22} \left[ C_{22} A_{22}^{-1} A_{21}\right]\\
&\quad + Q_{11} - A_{21}^* A_{22}^{-*} Q_{12}^* - Q_{12} A_{22}^{-1} A_{21} \\
&\quad + A_{21}^* A_{22}^{-*} \left[A_{22}^* \overline{P}_{22} + \overline{P}_{22} A_{22} + C_{12}^* \overline{P}_{11} C_{12} + C_{22}^* \overline{P}_{22} C_{22} + Q_{22}\right] A_{22}^{-1} A_{21}\\
&= A_s^* \overline{P}_{11} + \overline{P}_{11}A_s  + C_{1s}^* \overline{P}_{11} C_{1s} + C_{2s}^* \overline{P}_{22} C_{2s} + Q_s - \left[ \overline{F}_1 - \overline{F}_2 A_{22}^{-1}A_{21}\right]^* \overline{\Delta} \left[ \overline{F}_1 - \overline{F}_2 A_{22}^{-1}A_{21}\right].
\end{align*}
Now, observe that using \eqref{reduced6}, we have
\begin{align*}
\left[ I + \overline{F}_2 A_{22}^{-1} B_2\right]^* & \overline{\Delta} \left[ I + \overline{F}_2 A_{22}^{-1} B_2\right] \\
&= \overline{\Delta} + \overline{\Delta} \overline{F}_2 A_{22}^{-1} B_2 + B_2^* A_{22}^{-*} \overline{F}_2^* \overline{\Delta} + B_2^* A_{22}^{-*} \overline{F}_2^* \overline{\Delta} \overline{F}_2 A_{22}^{-1} B_2\\
&= \overline{\Delta} + \overline{\Delta} \overline{F}_2 A_{22}^{-1} B_2 + B_2^* A_{22}^{-*} \overline{F}_2^* \overline{\Delta} \\
&\quad + B_2^* A_{22}^{-*} \left[A_{22}^* \overline{P}_{22} + \overline{P}_{22} A_{22} + C_{12}^* \overline{P}_{11} C_{12} + C_{22}^* \overline{P}_{22} C_{22} + Q_{22}\right] A_{22}^{-1} B_2\\
&= R + D_1^* \overline{P}_{11} D_1 + D_2^* \overline{P}_{22} D_2 + B_2^* A_{22}^{-1} Q_{22} A_{22}^{-1} B_2 \\
&\quad + \left[ B_2^* \overline{P}_{22} + \overline{\Delta} \overline{F}_2 \right] A_{22}^{-1} B_2 + B_2^* A_{22}^{-*} \left[ \overline{F}_2^* \overline{\Delta} + \overline{P}_{22} B_2 \right] \\
&\quad + B_2^* A_{22}^{-*} \left[ C_{12}^* \overline{P}_{11} C_{12} + C_{22}^* \overline{P}_{22} C_{22} \right] A_{22}^{-1} B_2.
\end{align*}
Using \eqref{reducedoperator2} and completing the square, we have that
\begin{align*}
\left[ I + \overline{F}_2 A_{22}^{-1} B_2\right]^* & \overline{\Delta} \left[ I + \overline{F}_2 A_{22}^{-1} B_2\right] \\
&= R_s + D_1^* \overline{P}_{11} D_1 + D_2^* \overline{P}_{22} D_2  \\
&\quad - \left[ D_1^* \overline{P}_{11} C_{12} + D_2^* \overline{P}_{22} C_{22}  \right] A_{22}^{-1} B_2 - B_2^* A_{22}^{-*} \left[ D_1^* \overline{P}_{11} C_{12} + D_2^* \overline{P}_{22} C_{22}  \right]^* \\
&\quad + B_2^* A_{22}^{-*} \left[ C_{12}^* \overline{P}_{11} C_{12} + C_{22}^* \overline{P}_{22} C_{22} \right] A_{22}^{-1} B_2\\
&= R_s + D_{1s}^* \overline{P}_{11} D_{1s} + D_{2s}^* \overline{P}_{22} D_{2s}\\
&= \overline{\Delta}_s.
\end{align*}
This shows that $\overline{\Delta} > 0$ implies $\overline{\Delta}_s > 0$. Moreover, by Lemma \ref{lemma: invertible00}, we have that
\begin{equation}\label{temp91}
\overline{\Delta} = \left[ I + \overline{F}_2 A_{22}^{-1} B_2\right]^{-*} \overline{\Delta}_s\left[ I + \overline{F}_2 A_{22}^{-1} B_2\right]^{-1}.
\end{equation}
Applying \eqref{temp91} to $H$, we have that
\begin{align*}
H &= A_s^* \overline{P}_{11} + \overline{P}_{11}A_s  + C_{1s}^* \overline{P}_{11} C_{1s} + C_{2s}^* \overline{P}_{22} C_{2s} + Q_s \\
&\quad - \left[ \overline{F}_1 - \overline{F}_2 A_{22}^{-1}A_{21}\right]^* \left[ I + \overline{F}_2 A_{22}^{-1} B_2\right]^{-*} \overline{\Delta}_s\left[ I + \overline{F}_2 A_{22}^{-1} B_2\right]^{-1}\left[ \overline{F}_1 - \overline{F}_2 A_{22}^{-1}A_{21}\right].
\end{align*}
Lastly, we need to show that
\begin{equation}\label{temp94}
\overline{\Delta}_s\left[ I + \overline{F}_2 A_{22}^{-1} B_2\right]^{-1}\left[ \overline{F}_1 - \overline{F}_2 A_{22}^{-1}A_{21}\right]  = - M_s.
\end{equation}
To do so, we can use \eqref{temp91} to obtain
\begin{align*}
\overline{\Delta}_s\left[ I + \overline{F}_2 A_{22}^{-1} B_2\right]^{-1}\left[ \overline{F}_1 - \overline{F}_2 A_{22}^{-1}A_{21}\right] &= \left[ I + \overline{F}_2 A_{22}^{-1} B_2\right]^* \overline{\Delta} \left[ \overline{F}_1 - \overline{F}_2 A_{22}^{-1}A_{21}\right] \\
&= \overline{\Delta} \overline{F}_1 - \overline{\Delta} \overline{F}_2 A_{22}^{-1} A_{21} + B_2^* A_{22}^{-*} \overline{F}_2^* \overline{\Delta} \overline{F}_1\\
&\quad - B_2^* A_{22}^{-*} \overline{F}_2^* \overline{\Delta} \overline{F}_2 A_{22}^{-1} A_{21}.
\end{align*}
Applying \eqref{reducedoperator1}, \eqref{reducedoperator2}, \eqref{reduced5} and \eqref{reduced6}, we obtain
\begin{align*}
\overline{\Delta}_s  & \left[ I + \overline{F}_2 A_{22}^{-1} B_2\right]^{-1}\left[ \overline{F}_1 - \overline{F}_2 A_{22}^{-1}A_{21}\right] \\
&= - \left[ B_1^* \overline{P}_{11} + B_2^* \overline{P}_{12}^* + D_1^* \overline{P}_{11} C_{11} + D_2^* \overline{P}_{22} C_{21} \right] \\
&\quad +\left[ B_2^* \overline{P}_{22} + D_1^* \overline{P}_{11} C_{12} + D_2^* \overline{P}_{22} C_{22} \right] A_{22}^{-1} A_{21}\\
&\quad + B_2^* A_{22}^{-*} \left[ A_{12}^* \overline{P}_{11} + \overline{P}_{22} A_{21} + A_{22}^* \overline{P}_{12}^* + C_{12}^* \overline{P}_{11} C_{11} + C_{22}^* \overline{P}_{22} C_{21} + Q_{12}^*\right]\\
&\quad - B_2^* A_{22}^{-*} \left[ A_{22}^* \overline{P}_{22} + \overline{P}_{22} A_{22} + C_{12}^* \overline{P}_{11} C_{12} + C_{22}^* \overline{P}_{22} C_{22} + Q_{22}\right] A_{22}^{-1} A_{21}.
\end{align*}
Finally, by completing the squares, we obtain
\begin{align*}
\overline{\Delta}_s  \left[ I + \overline{F}_2 A_{22}^{-1} B_2\right]^{-1}\left[ \overline{F}_1 - \overline{F}_2 A_{22}^{-1}A_{21}\right] &= - B_s^* \overline{P}_{11} - D_{1s}^* \overline{P}_{11} C_{1s} - D_{2s}^* \overline{P}_{22}C_{2s} \\
&\quad - B_2^* A_{22}^{-*} Q_{22} A_{22}^{-1} A_{21} + B_2^* A_{22}^{-*} Q_{12}^*\\
&= - B_s^* \overline{P}_{11} - D_{1s}^* \overline{P}_{11} C_{1s} - D_{2s}^* \overline{P}_{22}C_{2s} \\
&\quad - B_2^* A_{22}^{-*} \left[ Q_{22} A_{22}^{-1} A_{21} - Q_{12}^* \right]\\
&= - B_s^* \overline{P}_{11} - D_{1s}^* \overline{P}_{11} C_{1s} - D_{2s}^* \overline{P}_{22}C_{2s} - L_s\\
&= -M_s.
\end{align*}
This completes the proof of Theorem \ref{theorem: equiv}-(1).

To show Theorem \ref{theorem: equiv}-(2), we can use a reversed argument to the above and use the invertibility of $A_{22} +B_2 \overline{F}_2$ to construct $\overline{P}_{12}$ as defined in \eqref{P12equiv}, which is clearly a solution to \eqref{reduced2}. 
\qed

\section{Well-posedness of the reduced differential-algebraic Riccati equation}\label{section: reducedDARE}

The goal of this section is to show the existence of a solution to the reduced DARE and consequently, the existence of a solution to the reduced system \eqref{reducedsystem} as well. We begin with a useful lemma on Lyapunov equations.

\begin{lemma}\label{lemma: difflyapunov}
Let $\Theta \in L^2([0,T]; \R^{k \times n_1})$. Denote $P \in C([0,T]; \mathbb{S}^{n_1})$ as the solution of the following Lyapunov equation 
\begin{equation}\label{difflyapunov}
\begin{split}
\begin{cases}
\frac{dP}{dt} + (A + B\Theta)^* P + P (A + B \Theta) + (C + D\Theta)^* P (C + D\Theta) \\
\qquad + \Theta^* R \Theta + S^*\Theta + \Theta^* S + Q = 0,\quad t\in [0,T],\\
P(T) = 0 \in \mathbb{S}^{n_1},
\end{cases}
\end{split}
\end{equation}
where the matrices $A, B, C, D, Q, R$ and $S$ are of appropriate dimensions. If the matrix
\begin{equation}\label{lyapunovconvex}
\mathcal{M}_1 := \begin{bmatrix}
Q & S^*\\
S & R
\end{bmatrix} > 0
\end{equation}
then, for every $\Theta \in L^2([0,T]; \R^{k \times n_1})$, the Lyapunov equation \eqref{difflyapunov} admits a unique solution $P \in C([0,T]; \mathbb{S}^{n_1})$ satisfying 
\begin{equation}
R + D^* P(t) D \geq \lambda I, \text{ and } P(t) \geq \alpha I,\quad \forall t\in [0,T],
\end{equation}
where $\alpha \geq 0$ and $\lambda > 0$.
\end{lemma}
\noindent \textit{Proof.} Consider the minimisation problem described by the quadratic cost functional
\begin{equation}
J(x;u) = \frac{1}{2} \E \int_0^T \left\langle \begin{bmatrix}
Q & S^*\\
S & R
\end{bmatrix} \begin{bmatrix}
X(t)\\ u(t)
\end{bmatrix},
\begin{bmatrix}
X(t)\\ u(t)
\end{bmatrix} \right\rangle dt, \quad u \in \mathbb{A}_T^k,
\end{equation}
subject to the linear stochastic differential equation
\begin{equation}
\begin{split}
\begin{cases}
dX(t) = [AX(t) + B u(t)]dt + [CX(t) + Du(t)] dW(t), \quad t\in [0,T],\\
X(0) = x.
\end{cases}
\end{split}
\end{equation}
From \eqref{lyapunovconvex}, we have that 
\begin{equation}
J(0;u) \geq \lambda \E \int_0^T |u(t)|^2 dt, \quad \forall u \in \mathbb{A}_T^k,
\end{equation}
where $\lambda > 0$. In addition, from Corollary 4.7 of \cite{freiling2004class} and Theorem 6.1 in Chapter 6 of \cite{yong1999stochastic}, we have that the value function defined on the time horizon $[s,T]$ satisfies 
\begin{equation}
V(s,x) \geq \alpha|x|^2,\quad \forall (s,x)\in [t,T]\times \R^{n_1},
\end{equation}
for some $\alpha \geq 0$. The result then follows from Proposition 4.4 of \cite{sun2016open}.
\qed

Next we derive some properties about the algebraic Riccati equation perturbed by the parameter $\overline{P}_{11}$.

\begin{lemma}\label{lemma: AREproperties}
Fix $\overline{P}_{11} \in \mathbb{S}^{n_1}_+$. Suppose that Assumption \ref{mainassumptions} holds. Then the algebraic Riccati equation
\begin{equation}\label{reducedAREfixed}
\begin{split}
\begin{cases}
A_{22}^* \overline{P}_{22} + \overline{P}_{22} A_{22} + C_{22}^* \overline{P}_{22} C_{22} + \left[Q_{22} + C_{12}^* \overline{P}_{11} C_{12} \right] \\
\quad- \left[ B_2^* \overline{P}_{22}  + D_{2}^* \overline{P}_{22} C_{22} + D_1^* \overline{P}_{11} C_{12} \right]^*\left[ \left(R + D_1^* \overline{P}_{11} D_1\right) + D_2^* \overline{P}_{22} D_2\right]^{-1}  \\
\qquad \left[ B_2^* \overline{P}_{22}  + D_{2}^* \overline{P}_{22} C_{22} + D_1^* \overline{P}_{11} C_{12} \right]= 0,\\
R + D_1^* \overline{P}_{11} D_1 + D_2^* \overline{P}_{22} D_2 > 0,
\end{cases}
\end{split}
\end{equation}
admits a solution $\overline{P}_{22}^+ \in \mathbb{S}^{n_2}_{++}$ unique in the set of positive semi-definite solutions. In addition:
\begin{enumerate}
\item The eigenvalues of the operator $A_{22} + B_2 \overline{F}_2$ have negative real parts where
\begin{equation}\label{temp20}
\begin{split}
\overline{F}_2(\overline{P}_{11}) &= - (R + D_1^* \overline{P}_{11} D_1 + D_2^* \overline{P}_{22}^+ D_2)^{-1} \\
&\qquad \qquad\left[B_2^* \overline{P}_{22}^+ + D_2^*\overline{P}_{22}^+ C_{22}  + D_1^* \overline{P}_{11} C_{12}  \right].
\end{split}
\end{equation}
This implies that there exists constants $M \geq 1$ $\gamma > 0$ such that
\begin{equation}\label{expstabG}
|e^{t( A_{22} + B_2 \overline{F}_2)}| \leq M e^{-\gamma t},\quad t\geq 0.
\end{equation}
\item We can write $\overline{P}_{22}^+ = h_2(\overline{P}_{11})$ where $h_2: \mathbb{S}^{n_1} \rightarrow \mathbb{S}^{n_2}$ is an increasing continuously differentiable function. That is, $\overline{P}_{11} \leq \overline{P}_{11}'$ implies $h_2(\overline{P}_{11})\leq h_2(\overline{P}_{11}')$. Moverover, $h_2$ bounded by the linear function $h_2':\mathbb{S}^{n_1} \rightarrow \mathbb{S}^{n_2}$ defined as
\begin{equation}\label{alglyapunovsoln}
h_2'(\overline{P}_{11}) := \E \int_0^\infty \Phi(t)^* \left[Q_{22} + C_{12}^* \overline{P}_{11} C_{12} \right] \Phi(t) dt
\end{equation}
where $\Phi$ is the continuous $\mathcal{F}$-adapted solution to the SDE
\begin{equation}
\begin{split}
\begin{cases}
d\Phi(t) = A_{22} \Phi(t) dt + C_{22} \Phi(t) dW(t), t\geq 0,\\
\Phi(0) = I,
\end{cases}
\end{split}
\end{equation}
satisfying
\begin{align*}
\E \int_0^\infty |\Phi(t)|^2 dt < \infty.
\end{align*}
By bounded, we mean that $h_2(\overline{P}_{11}) \leq h_2'(\overline{P}_{11})$ for all $\overline{P}_{11} \in \mathbb{S}^{n_1}_+$.
\item The differential Riccati equation with an initial condition
\begin{equation}\label{reducedAREfixeddiff}
\begin{split}
\begin{cases}
\frac{dP_{22}(t)}{dt} = A_{22}^* P_{22}(t) + P_{22}(t) A_{22} + C_{22}^* P_{22}(t) C_{22} + \left[Q_{22} + C_{12}^* \overline{P}_{11} C_{12} \right] \\
\quad- \left[ B_2^* P_{22}(t)  + D_{2}^* P_{22}(t) C_{22} + D_1^* \overline{P}_{11} C_{12} \right]^*\left[R + D_1^* \overline{P}_{11} D_1 + D_2^* P_{22}(t) D_2\right]^{-1}  \\
\qquad \left[ B_2^* P_{22} (t) + D_{2}^* P_{22}(t) C_{22} + D_1^* \overline{P}_{11} C_{12} \right],\\
P_{22}(0) \geq 0,\\
R + D_1^* \overline{P}_{11} D_1 + D_2^* P_{22}(t) D_2 > 0, \quad t\in [0,T]\ a.e.,
\end{cases}
\end{split}
\end{equation}
admits a unique solution $P_{22} \in C([0,\infty);\mathbb{S}^{n_2}_+)$ and moreover satisfies $\lim_{t\rightarrow \infty} P_{22}(t) = \overline{P}_{22}^+$. In addition, the equilibrium $\overline{P}_{22}^+$ is exponentially stable. That is, for arbitrary $p> 0$,
\begin{equation}\label{expstab0}
| P_{22}(t) - \overline{P}_{22}^+ | \leq M_P e^{- \frac{3\gamma}{2} t} |P_{22}(0) - \overline{P}_{22}^+|,\quad \forall |\overline{P}_{11}| \leq p,\ \forall t \geq 0,
\end{equation}
for some positive constant $M_P$, which may depend on $p$, and $\gamma$ is as defined in \eqref{expstabG}.
\end{enumerate}
\end{lemma}
\noindent \textit{Proof.} From Assumption \ref{mainassumptions}, the matrices $Q_{11}$ and $R_s$ are positive definite. This implies that the matrix
\begin{equation}\label{temp505}
\begin{split}
\mathcal{M}_2(\overline{P}_{11}) := \begin{bmatrix}
Q_{22} + C_{12}^* \overline{P}_{11} C_{12} & C_{12}^* \overline{P}_{11} D_1\\
D_1^* \overline{P}_{11} C_{12} & R_s + D_1^* \overline{P}_{11} D_1
\end{bmatrix} &= \begin{bmatrix}
Q_{22} & 0\\
0 & R_s
\end{bmatrix} + \begin{bmatrix}
C_{12}^*\\
D_1^*
\end{bmatrix} \overline{P}_{11}
\begin{bmatrix}
C_{12}^*\\
D_1^*
\end{bmatrix}^*\\
&\geq 
\begin{bmatrix}
Q_{22} & 0\\
0 & R_s
\end{bmatrix} >0.
\end{split}
\end{equation}
Hence, by Theorem 1 of Section II.E of \cite{rami2000linear} and Theorem 3.3, Corollary 5.3, Lemma 5.5 and Lemma 5.8 of \cite{freiling2004class}, the reduced algebraic Riccati equation \eqref{reducedAREfixed} admits a unique positive solution $\overline{P}_{22}^+$ such that the eigenvalues of $A_{22} + B_2 \overline{F}_2$ have negative real parts and $[A_{22} + B_2 \overline{F}_2, C_{22} + D_2 \overline{F}_2]$ is $L^2$-stable. Moreover, it is unique within the class of positive semi-definite solutions.   

Next, we show that the solution $\overline{P}_{22}^+$ can be written as a differentiable function $h_2(\overline{P}_{11})$. To this end, we will apply the Implicit Function Theorem. For convenience, we rewrite \eqref{reducedAREfixed} in terms of the solution $(\overline{P}^+_{22},\overline{F}_2)$
\begin{equation}\label{reducedAREfixed2}
\begin{split}
\begin{cases}
A_{22}^* \overline{P}_{22}^+ + \overline{P}_{22}^+ A_{22} + C_{22}^* P_{22}^+ C_{22} - \overline{F}_2^* \Delta_2 \overline{F}_2 + \left[Q_{22} + C_{12}^* \overline{P}_{11} C_{12} \right] = 0,\\
B_2^* \overline{P}_{22}^+ + D_1^* \overline{P}_{11} C_{12} + D_2^* \overline{P}_{22}^+ C_{22} + \Delta_2 \overline{F}_2 = 0.
\end{cases}
\end{split}
\end{equation}
with 
\begin{equation}\label{reducedAREfixedcon}
\Delta_2 = R + D_1^* \overline{P}_{11} D_1 + D_2^* \overline{P}_{22}^+ D_2 > 0.
\end{equation} 
Let us denote the system \eqref{reducedAREfixed2} in the compact form 
\begin{equation}
\mathbf{F}(\overline{P}_{22}^+, \overline{F}_2; \overline{P}_{11}) = 0.
\end{equation}
We observe that $\mathbf{F}: \mathbb{S}^{n_2} \times \R^{k \times n_2} \times \mathbb{S}^{n_1} \rightarrow  \mathbb{S}^{n_2} \times \R^{k \times n_2}$ is a differentiable function. We need to check that the mapping 
\begin{equation}
\begin{split}
(Z_1,Z_2) \rightarrow &\frac{\partial \mathbf{F}(\overline{P}_{22}^+, \overline{F}_2; \overline{P}_{11})}{\partial (\overline{P}_{22}^+, \overline{F}_2)} (Z_1,Z_2) \\
&= \lim_{h \rightarrow 0} \frac{1}{h} \left[ \mathbf{F}(\overline{P}_{22}^+ + h Z_1, \overline{F}_2 + h Z_2; \overline{P}_{11}) - \mathbf{F}(\overline{P}_{22}^+, \overline{F}_2; \overline{P}_{11} ) \right]
\end{split}
\end{equation}
is an isomorphism. As the above differential mapping takes values from the finite-dimensional vector space $\mathbb{S}^{n_2} \times \R^{k \times n_2}$ to itself, it is enough to check that this mapping is injective. That is, the unique solution to 
\begin{equation}\label{injectiveeqn}
\frac{\partial \mathbf{F}(\overline{P}_{22}^+, \overline{F}_2; \overline{P}_{11})}{\partial (\overline{P}_{22}^+, \overline{F}_2)} (Z_1,Z_2) = 0
\end{equation}
is $(Z_1,Z_2) = 0$. Evaluating \eqref{injectiveeqn} gives the equations
\begin{equation}\label{temp1}
A_{22}^* Z_1 + Z_1 A_{22} + C_{22}^* Z_1 C_{22} - (\overline{F}_2)^* D_2^* Z_1 D_2 \overline{F}_2 - (\overline{F}_2)^* \Delta_2 Z_2 - Z_2^* \Delta_2 \overline{F}_2 = 0
\end{equation} 
\begin{equation}\label{temp2}
B_2^* Z_1 + D_2^* Z_1 C_{22} + \Delta_2 Z_2 + D_2^* Z_1 D_2 \overline{F}_2 = 0.
\end{equation}
Using \eqref{temp2} we can eliminate the variable $Z_2$ from \eqref{temp1} yields the Lyapunov equation
\begin{equation}\label{lyapunov}
\left[ A_{22} + B_2 \overline{F}_2 \right]^* Z_1 + Z_1 \left[ A_{22} + B_2 \overline{F}_2 \right] + \left[ C_{22} + D_2 \overline{F}_2 \right]^* Z_1 \left[ C_{22} + D_2 \overline{F}_2 \right] = 0.
\end{equation}
It is obvious that $Z_1 = 0$ is a solution to \eqref{lyapunov}. Since the eigenvalues of $ A_{22} + B_2 \overline{F}_2$ have negative real parts and $[A_{22} + B_2 \overline{F}_2, C_{22} + D_2 \overline{F}_2]$ is $L^2$-stable, applying Lemma 2.2 of \cite{sun2018stochastic} gives the uniqueness of $Z_1 = 0$. Consequently, $Z_2 = 0$ must also be the unique solution to \eqref{temp2}. Hence, by the Implicit Function Theorem (see \cite{edwards2012advanced}), we can write $\overline{P}_{22}^+ = h_2(\overline{P}_{11})$ for some differentiable function $h_2$. To see that $h_2$ is also an increasing function (in the definiteness sense), let us assume that $\widetilde{P}_{11} \geq \overline{P}_{11}$. It is clear that $\mathcal{M}_2(\widetilde{P}_{11}) \geq \mathcal{M}_2 (\overline{P}_{11})$. Hence, by Theorem 5.3 of \cite{freiling2001basic},  $h_2(\widetilde{P}_{11}) \geq h_2 (\overline{P}_{11}).$ 

We will now show that $h_2$ is bounded from above by a linear function $h_2'$. That is, $h_2(\overline{P}_{11}) \leq h_2'(\overline{P}_{11})$ for all $\overline{P}_{11} \in \mathbb{S}^{n_1}_+$. For fixed $\overline{P}_{11} \in \mathbb{S}^{n_1}_+$, consider the Lyapunov equation
\begin{equation}\label{alglyapunov}
A_{22}^* \overline{P}'_{22} + \overline{P}'_{22} A_{22} + C_{22}^* \overline{P}'_{22} C_{22} + \left[Q_{22} + C_{12}^* \overline{P}_{11} C_{12} \right] = 0.
\end{equation}
Since $[A_{22}, C_{22}]$ is $L^2$-stable, we can apply Lemma 2.2 of \cite{sun2018stochastic}, to show that the solution to the Lyapunov equation \eqref{alglyapunov} is given by
\begin{equation}\label{alglyapunovsolnwrong}
\overline{P}_{22}' = h_2'(\overline{P}_{11}) := \E \int_0^\infty \Phi(t)^* \left[Q_{22} + C_{12}^* \overline{P}_{11} C_{12} \right] \Phi(t) dt
\end{equation}
where $\Phi$ is the continuous $\mathcal{F}$-adapted solution to the SDE
\begin{equation}
\begin{split}
\begin{cases}
d\Phi(t) = A_{22} \Phi(t) dt + C_{22} \Phi(t) dW(t), t\geq 0,\\
\Phi(0) = I,
\end{cases}
\end{split}
\end{equation}
satisfying
\begin{align*}
\E \int_0^\infty |\Phi(t)|^2 dt < \infty.
\end{align*}
Clearly $h'$ is linear in $\overline{P}_{11}$. Now, let $Z = \overline{P}_{22} - \overline{P}_{22}'$. From \eqref{reducedAREfixed} and \eqref{alglyapunov}, we have that
\begin{equation}
A_{22}^* Z + Z A_{22} + C_{22}^* Z C_{22} > 0.
\end{equation}
Again, since $[A_{22}, C_{22}]$ is $L^2$-stable, we can apply Lemma 2.2 of \cite{sun2018stochastic} to show that $Z < 0$. Hence for every $\overline{P}_{11}\in  \mathbb{S}^{n_1}$, the function $h(\overline{P}_{11})$ is bounded from above by the linear function $h_2'(\overline{P}_{11})$.

Now consider the differential Riccati equation \eqref{reducedAREfixeddiff}. Since \eqref{temp505} holds, the existence and uniqueness of a solution $P_{22} \in C([0,T]; \mathbb{S}^{n_2}_+)$ follows from Theorem 7.2 in Chapter 6 of \cite{yong1999stochastic}. The convergence result follows from Theorem 6.3 of \cite{freiling2004class}. To show the exponential stability of the equilibrium $\overline{P}_{22}^+$, let us set $Y(t) = P_{22}(t) - \overline{P}_{22}^+$ and
\begin{equation}\label{temp21}
\begin{split}
F_2(\overline{P}_{11})(t) &= - \left[R + D_1^* \overline{P}_{11} D_1 + D_2^* P_{22}(t) D_2\right]^{-1}  \\
&\qquad \quad \left[ B_2^* P_{22} (t) + D_{2}^* P_{22}(t) C_{22} + D_1^* \overline{P}_{11} C_{12} \right],\quad t \geq 0.
\end{split}
\end{equation}
Observe that
\begin{equation}\label{temp22}
|F_2(\overline{P}_{11})(t) - \overline{F}_2(\overline{P}_{11})| \leq \overline{M} |P_{22}(t) - \overline{P}_{22}^+|,\quad \forall t\geq 0,
\end{equation}
for some positive constant $\overline{M}$ independent of $t$. Using a similar parametrisation as \eqref{reducedAREfixed2}, we can write $Y(t)$ as the solution to the differential equation
\begin{equation}
\begin{split}
\begin{cases}
\frac{dY}{dt} = A_{22}^* Y + Y A_{22} + C_{22}^* Y C_{22} - F_2^* (R + D_1^* \overline{P}_{11} D_1 + D_2^* P_{22} D_2) F_2\\
\qquad \qquad  + \overline{F}_2^* (R + D_1^* \overline{P}_{11} D_1 + D_2^* \overline{P}_{22}^+ D_2) \overline{F}_2,\\
Y(0) = P_{22}(0) - \overline{P}_{22}^+.
\end{cases}
\end{split}
\end{equation}
Adding and subtracting terms, we have that
\begin{align*}
\frac{dY}{dt} 
&=  \left[A_{22}+ B_2 \overline{F}_2\right]^* Y + Y \left[A_{22}+ B_2 \overline{F}_2\right] + C_{22}^* Y C_{22}  - \overline{F}_2^* B_2^* Y - Y B_2 \overline{F}_2\\
&\quad  - F_2^* (R + D_1^* \overline{P}_{11} D_1 + D_2^* P_{22} D_2) F_2+ \overline{F}_2^* (R + D_1^* \overline{P}_{11} D_1 + D_2^* \overline{P}_{22}^+ D_2) \overline{F}_2.
\end{align*}
From \eqref{temp20} and \eqref{temp21}, we can simplify the above to
\begin{align*}
\frac{dY}{dt} 
&=  \left[A_{22}+ B_2 \overline{F}_2\right]^* Y + Y \left[A_{22}+ B_2 \overline{F}_2\right] + C_{22}^* Y C_{22}  - Y B_2 \overline{F}_2\\
&\quad  + F_2^* \left[ D_1^* \overline{P}_{11} C_{12} + D_2^* P_{22} C_{22} \right] - \overline{F}_2^* \left[ D_1^* \overline{P}_{11} C_{12} + D_2^* \overline{P}_{22}^+ C_{22} \right]\\
&=  \left[A_{22}+ B_2 \overline{F}_2\right]^* Y + Y \left[A_{22}+ B_2 \overline{F}_2\right] + C_{22}^* Y C_{22}  - Y B_2 \overline{F}_2\\
&\quad  + F_2^*  D_2^* Y C_{22}  + \left[F_2 - \overline{F}_2\right]^* \left[ D_1^* \overline{P}_{11} C_{12} + D_2^* \overline{P}_{22}^+ C_{22} \right].
\end{align*}
Applying a variation of constants formula,
we obtain
\begin{align*}
Y(t) &=  e^{(A_{22} + B_2 \overline{F}_2)^* t} \left[P_{22}(0) - \overline{P}_{22}^+\right] e^{(A_{22} + B_2 \overline{F}_2) t}\\
&\quad + \int_0^t e^{(A_{22} + B_2 \overline{F}_2)^* (t-s)} \left[ C_{22}^* Y(s) C_{22}  - Y(s) B_2 \overline{F}_2 + F_2^*  D_2^* Y(s) C_{22}  \right] e^{(A_{22} + B_2 \overline{F}_2) (t-s)} ds\\
&\quad + \int_0^t e^{(A_{22} + B_2 \overline{F}_2)^* (t-s)}  \left[F_2(s) - \overline{F}_2\right]^* \left[ D_1^* \overline{P}_{11} C_{12} + D_2^* \overline{P}_{22}^+ C_{22} \right]e^{(A_{22} + B_2 \overline{F}_2) (t-s)} ds.
\end{align*}
Hence applying \eqref{temp22} and the fact that the eigenvalues of $A_{22} + B_2 \overline{F}_2$ have negative real parts, we obtain 
\begin{align*}
|Y(t)| &\leq  K_1 e^{- 2 \gamma t} |P_{22}(0) - \overline{P}_{22}^+| + K_2 \int_0^t e^{-2\gamma (t-s)} |Y(s)| ds
\end{align*}
for some positive constants $K_1, K_2$ and $\gamma$, which may depend on $p$. So we have that,
\begin{align*}
e^{\frac{3\gamma}{2} t} |Y(t)| &\leq  K_1 |P_{22}(0) - \overline{P}_{22}^+| + K_2 \int_0^t e^{-\frac{\gamma}{2}(t-s)} \left[e^{\frac{3\gamma}{2} s}  |Y(s)| \right] ds
\end{align*}
Applying Gronwall's inequality (see Theorem 15 of \cite{dragomir2003some}), we obtain
\begin{align*}
|Y(t)| \leq K_1 e^{\frac{2K_2}{\gamma}} e^{- \frac{3\gamma}{2} t}|P_{22}(0) - \overline{P}_{22}^+|
\end{align*}
as required.
\qed

\begin{theorem}\label{theorem: DARE}
Suppose that Assumption \ref{mainassumptions} holds. Then the reduced differential-algebraic Riccati equation admits a unique solution $(\overline{P}_{11}^+, \overline{P}_{22}^+) \in C([0,T]; \mathbb{S}^{n_1}_+) \times C([0,T]; S^{n_2}_{++})$ such that $\overline{P}_{22}^+$ is unique amongst the set of positive semi-definite solutions. Moreover, we can write $\overline{P}_{22}^+(t) = h_2(\overline{P}_{11}^+(t)), t\in [0,T]$ where $h_2$ is the continuously differentiable function defined in Lemma \ref{lemma: AREproperties}. In addition, the eigenvalues of $A_{22} + B_2 \overline{F}_2^+(t)$ have negative real parts for all $t\in [0,T]$ where
\begin{equation}
\begin{split}
\overline{F}_2^+(t) &= - (R + D_1^* \overline{P}_{11}^+(t) D_1 + D_2^* \overline{P}_{22}^+(t) D_2)^{-1} \\
&\qquad \qquad\left[B_2^* \overline{P}_{22}^+(t) + D_1^* \overline{P}_{11}^+(t) C_{12} + D_2^* \overline{P}_{22}^+(t) C_{22}\right], \quad t\in [0,T].
\end{split}
\end{equation}
\end{theorem}
\noindent\textit{Proof.} From Lemma \ref{lemma: AREproperties}, we have shown that for fixed $\overline{P}_{11}$, the reduced algebraic Riccati equation \eqref{reducedARE} admits a unique solution $\overline{P}_{22} = h_2(\overline{P}_{11})$ such that $h_2$ is an increasing and continuously differentiable function bounded by a linear function $h_2'$. Moreover, the eigenvalues of the operator $A_{22} + B_2 \overline{F}_2$ have negative real parts where
\begin{equation}
\begin{split}
\overline{F}_2(\overline{P}_{11}) &= - (R + D_1^* \overline{P}_{11} D_1 + D_2^* h_2(\overline{P}_{11}) D_2)^{-1} \\
&\qquad \qquad\left[B_2^* h_2(\overline{P}_{11}) + D_1^* \overline{P}_{11} C_{12} + D_2^* h_2(\overline{P}_{11}) C_{22}\right].
\end{split}
\end{equation}

We now turn to the reduced differential Riccati equation \eqref{reducedDRE} with the condition \eqref{convexDRE} and replace each instance of $\overline{P}_{22}$ with $h_2(\overline{P}_{11})$. That is,
\begin{subnumcases}{}
\label{reducedDREfixed}
\begin{split}
&\frac{d\overline{P}_{11}}{dt} + A_s^* \overline{P}_{11} + \overline{P}_{11} A_2 + C_{1s}^* \overline{P}_{11} C_{1s} + C_{2s}^* h_2(\overline{P}_{11}) C_{2s} + Q_s\\
&\quad \qquad- M_s^* (R_s + D_{1s}^* \overline{P}_{11} D_{1s} + D_{2s}^* h_2(\overline{P}_{11})  D_{2s})^{-1} M_s = 0, \quad \overline{P}_{11}(T) = 0,
\end{split}\\
\label{convexDREfixed}
R_s + D_{1s}^* \overline{P}_{11}(t) D_{1s} + D_{2s}^* h_2(\overline{P}_{11}(t))  D_{2s} > 0, \quad \forall t\in [0,T]\ a.e.,
\end{subnumcases}
where $M_s = B_s^* \overline{P}_{11} + D_{1s}^* \overline{P}_{11} C_{1s} + D_{2s}^* h_2(\overline{P}_{11}) C_{2s} + L_s$. We will prove the existence of a solution to \eqref{reducedDREfixed}-\eqref{convexDREfixed} by an iterative scheme inspired by Theorem 4.5 of \cite{sun2016open}. 

Consider the differential Lyapunov equation
\begin{equation}\label{lyapunov0}
\begin{split}
\begin{cases}
\frac{d\overline{P}_{11}'}{dt} + A_s \overline{P}_{11}' + \overline{P}_{11}' A_s + C_{1s}^* \overline{P}_{11}' C_{1s} + Q_s = 0\\
P'(T) = 0.
\end{cases}
\end{split}
\end{equation}
We see that
\begin{align*}
\mathcal{M}_1^0 :&= \begin{bmatrix}
Q_s & L_s^*\\
L_s & R_s
\end{bmatrix} \\
&= \begin{bmatrix}
0 & 0\\
0 & R
\end{bmatrix} + \begin{bmatrix}
Q_s & L_s^*\\
L_s & B_2^* A_{22}^{-*} Q_{22} A_{22}^{-1} B_2
\end{bmatrix}\\
&= \begin{bmatrix}
0 & 0\\
0 & R
\end{bmatrix} + 
\begin{bmatrix}
I & 0\\
0 & B_2
\end{bmatrix}^* \begin{bmatrix}
Q_s & \left( A_{21}^* A_{22}^{-*} Q_{22} - Q_{12} \right) A_{22}^{-1}\\
A_{22}^{-*} \left( Q_{22} A_{22}^{-1} A_{21}  - Q_{12}^* \right)  & A_{22}^{-*} Q_{22} A_{22}^{-1} 
\end{bmatrix} \begin{bmatrix}
I & 0\\
0 & B_2
\end{bmatrix}.
\end{align*}
By Assumption \ref{mainassumptions}-(1) and Lemma 2.3 of \cite{freiling2004class}, the matrix 
\[
\begin{bmatrix}
Q_s & \left( A_{21}^* A_{22}^{-*} Q_{22} - Q_{12} \right) A_{22}^{-1}\\
A_{22}^{-*} \left( Q_{22} A_{22}^{-1} A_{21}  - Q_{12}^* \right)  & A_{22}^{-*} Q_{22} A_{22}^{-1} 
\end{bmatrix}
\]
is positive definite. Hence $\mathcal{M}_1^0$ is too. Thus, applying Lemma \ref{lemma: difflyapunov} with $\Theta = 0$, we have that \eqref{lyapunov0} admits a unique solution $\overline{P}'_{11} \in C([0,T]; \mathbb{S}^{n_1})$ satisfying
\begin{equation}
R_s + D_{1s}^* \overline{P}_{11}'(t) D_{1s} \geq \lambda I,\quad \overline{P}_{11}'(t) \geq \alpha_0 I,\quad \forall t\in [0,T]\ a.e.,
\end{equation}
for some $\alpha_0 \geq 0$ and $\lambda > 0$. Moreover, as $h_2$ is positive definite, we have that
\begin{equation}\label{lyapunovinq0.5}
R_s + D_{1s}^* \overline{P}_{11}'(t) D_{1s} +  D_{2s}^* h_2(\overline{P}_{11}'(t))  D_{2s} \geq \lambda I,\quad \forall t\in [0,T] \ a.e.
\end{equation}

Now consider the perturbed differential Lyapunov equation
\begin{equation}\label{lyapunov0.5}
\begin{split}
\begin{cases}
\frac{d\overline{P}_{11}^0}{dt} + A_s \overline{P}_{11}^0 + \overline{P}_{11}^0 A_s + C_{1s}^* \overline{P}_{11}^0 C_{1s} + Q_s + C_{2s}^* h_2'(\overline{P}_{11}^0) C_{2s} = 0, \\
P^0(T) = 0,
\end{cases}
\end{split}
\end{equation}
where $h'_2$ is the positive definite linear function defined in \eqref{alglyapunovsoln}. By Theorem 2.1 of \cite{wonham1968matrix}, \eqref{lyapunov0.5} admits a unique solution $\overline{P}_{11}^0 \geq 0$. Moreover, by taking the difference between \eqref{lyapunov0} and \eqref{lyapunov0.5} and using the positiveness of $h'$, Theorem 2.1 of \cite{wonham1968matrix} implies that $\overline{P}_{11}^0 \geq \overline{P}_{11}'$. Hence, by the increasing property of $h_2$ from Lemma \ref{lemma: AREproperties}, the inequality \eqref{lyapunovinq0.5} implies that
\begin{equation}\label{lyapunovinq0.56}
R_s + D_{1s}^* \overline{P}_{11}^0(t) D_{1s} +  D_{2s}^* h_2(\overline{P}_{11}^0(t))  D_{2s} \geq \lambda I,\quad \forall t\in [0,T]\ a.e.
\end{equation}

We now proceed via induction: For $i = 0,1,2,\ldots$, define 
\begin{equation}\label{lyapunov1}
\begin{split}
\begin{cases}
\Gamma_0 &= 0,\quad \Gamma_i = h_2(\overline{P}_{11}^{i-1}) - h_2(\overline{P}_{11}^{i}), \quad i \geq 1,\\
\Theta_i &= - \left(R_s + D_{1s}^* \overline{P}_{11}^i D_{1s} +  D_{2s}^* \left[h_2(\overline{P}_{11}^i) + \Gamma_i\right]  D_{2s}\right)^{-1}\\
&\qquad \qquad \left[B_s^* \overline{P}_{11}^i + D_{1s}^* \overline{P}_{11}^i C_{1s} + D_{2s}^* \left[h_2(\overline{P}_{11}^i) + \Gamma_i\right] C_{2s} + L_s \right],\\
A_i &= A_s + B_s \Theta_i,\quad C_i = C_{1s} + D_{1s} \Theta_i,\\
Q_i &= Q_s + C_{2s}^* \left[h_2(\overline{P}_{11}^i) + \Gamma_i\right] C_{2s},\quad R_i = R_s + D_{2s}^* \left[h_2(\overline{P}_{11}^i) + \Gamma_i\right]  D_{2s},\\
S_i &= L_s + D_{2s}^* \left[h_2(\overline{P}_{11}^i) + \Gamma_i\right] C_{2s}.
\end{cases}
\end{split}
\end{equation}
Let $\overline{P}_{11}^{i+1}$, given $\overline{P}_{11}^i$, be the solution to the Lyapunov equation 
\begin{equation}
\begin{split}
\begin{cases}
\frac{d\overline{P}_{11}^{i+1}}{dt} + A_i^* \overline{P}_{11}^{i+1}+ \overline{P}_{11}^{i+1} A_i + C_i^* \overline{P}_{11}^{i+1} C_i + \Theta_i^* R_i \Theta_i + S_i^*\Theta_i + \Theta_i^* S_i + Q_i = 0,\quad t\in [0,T],\\
\overline{P}_{11}^{i+1}(T) = 0.
\end{cases}
\end{split}
\end{equation}
Notice that
\begin{align*}
h_2(\overline{P}_{11}^i(t)) + \Gamma_i(t) =
\begin{cases}
 h_2(\overline{P}_{11}^0(t)) ,\quad i = 0,\\
 h_2(\overline{P}_{11}^{i-1}(t)),\quad i = 1,2,3,\ldots.
 \end{cases}
\end{align*}
Observe that the matrix
\begin{equation}
\mathcal{M}_1^i = \begin{bmatrix}
Q_i & S_i^*\\
S_i & R_i
\end{bmatrix} = \mathcal{M}_1^0 + \begin{bmatrix}
C_{2s}^* \\
D_{2s}^*
\end{bmatrix} \left[h_2(\overline{P}_{11}^i(t)) + \Gamma_i(t)\right]\begin{bmatrix}
C_{2s}^* \\
D_{2s}^*
\end{bmatrix}^*
\end{equation}
is positive definite. Hence, for $i = 1,2,3,\ldots$, given $\overline{P}_{11}^i$, we can apply Lemma \ref{lemma: difflyapunov} to ensure that the Lyapunov equation \eqref{lyapunov1} admits a unique solution $\overline{P}_{11}^{i+1} \in C([0,T]; \mathbb{S}^{n_1})$ such that
\begin{equation}\label{convexiteration}
R_s + D_{1s}^* \overline{P}_{11}^{i+1}(t) D_{1s} +  D_{2s}^* h_2(\overline{P}_{11}^{i-1}(t))  D_{2s} \geq \lambda_i I,\quad\text{and}\quad \overline{P}_{11}^{i+1}(t) \geq \alpha_i I,\quad \forall t\in [0,T]\ a.e.,
\end{equation}
for some $\lambda_i > 0$ and $\alpha_i \geq 0$. However as $\mathcal{M}_1^i > \mathcal{M}_1^0$, we observe from the proof of Lemma \ref{lemma: difflyapunov} that $\lambda_i > \lambda > 0$.

Next we show that the sequence $\{ \overline{P}_{11}^{i}\}_{i = 1}^\infty$ is uniformly bounded in $C([0,T]; \mathbb{S}^{n_1})$. To do so, for $i = 0,1,2,\ldots$, set
\begin{equation}
K_i = \overline{P}_{11}^{i} - \overline{P}_{11}^{i+1},\quad \text{and} \quad \Lambda_i = \Theta_{i-1} - \Theta_i.
\end{equation}
Hence we have that $K_i(T) = 0$ and for $i = 1,2,3,\ldots$,
\begin{equation}\label{temp3}
\begin{split}
-\frac{dK_i}{dt} &=   A_{i-1}^* \overline{P}_{11}^{i}+ \overline{P}_{11}^{i} A_{i-1} + C_{i-1}^* \overline{P}_{11}^{i} C_{i-1} \\
&\quad + \Theta_{i-1}^* R_{i-1} \Theta_{i-1} + S_{i-1}^*\Theta_{i-1} + \Theta_{i-1}^* S_{i-1} + Q_{i-1}\\
&\quad-A_i^* \overline{P}_{11}^{i+1}- \overline{P}_{11}^{i+1} A_i -C_i^* \overline{P}_{11}^{i+1} C_i  - \Theta_i^* R_i \Theta_i - S_i^*\Theta_i - \Theta_i^* S_i - Q_i\\
&= A_i^* K_i + K_i A_i + C_i^* K_i C_i + (A_{i-1} - A_i)^* \overline{P}_{11}^{i} + \overline{P}_{11}^i (A_{i - 1} - A_i)\\
&\quad +C_{i-1}^* \overline{P}_{11}^{i} C_{i-1}  - C_{i}^* \overline{P}_{11}^{i} C_{i}  + \Theta_{i-1}^* R_{i-1} \Theta_{i-1} - \Theta_i^* R_i \Theta_i \\
&\quad + S_{i-1}^*\Theta_{i-1} + \Theta_{i-1}^* S_{i-1} - S_i^*\Theta_i - \Theta_i^* S_i + Q_{i-1}- Q_i.
\end{split}
\end{equation}
By some algebraic manipulation, we can show that
\begin{equation}
\begin{split}
\begin{cases}
A_{i-1} - A_i = B_s \Lambda_i, \quad C_{i-1} -C_i = D_{1s} \Lambda_i,\\
C_{i-1}^* \overline{P}_{11}^{i} C_{i-1}  - C_{i}^* \overline{P}_{11}^{i} C_{i} = \Lambda_i^* D_{1s}^* \overline{P}_{11}^i D_{1s} \Lambda_i + C_i^* \overline{P}_{11}^i D_{1s} \Lambda_i + \Lambda_i^* D^*_{1s} \overline{P}_{11}^i C_i,\\
\Theta_{i-1}^* R_{i-1} \Theta_{i-1} - \Theta_i^* R_i \Theta_i = \Theta_{i-1}^* D_{2s}^* \Gamma_{i-1} D_{2s} \Theta_{i-1} + \Lambda_i^* R_i \Lambda_i + \Lambda_i^* R_i \Theta_i + \Theta_i^* R_i \Lambda_i,\\
S_{i-1}^*\Theta_{i-1} - S_i^*\Theta_i = C_{2s}^* \Gamma_{i-1} D_{2s} \Theta_{i-1} + S_i^* \Lambda_i,\\
Q_{i-1} - Q_i = C_{2s}^* \Gamma_{i-1} C_{2s}.
\end{cases}
\end{split}
\end{equation}
Applying these expressions to \eqref{temp3}, we have that
\begin{equation}\label{temp4}
\begin{split}
-\frac{dK_i}{dt}&= A_i^* K_i + K_i A_i + C_i^* K_i C_i + \Lambda_i^* B_s^* \overline{P}_{11}^{i} + \overline{P}_{11}^i  B_s \Lambda_i\\
&\quad +\Lambda_i^* D_{1s} \overline{P}_{11}^i D_{1s} \Lambda_i + C_i^* \overline{P}_{11}^i D_{1s} \Lambda_i + \Lambda_i^* D^*_{1s} \overline{P}_{11}^i C_i\\
&\quad  + \Theta_{i-1}^* D_{2s}^* \Gamma_{i-1} D_{2s} \Theta_{i-1} + \Lambda_i^* R_i \Lambda_i + \Lambda_i^* R_i \Theta_i + \Theta_i^* R_i \Lambda_i \\
&\quad + C_{2s}^* \Gamma_{i-1} D_{2s} \Theta_{i-1} + S_i^* \Lambda_i + \Theta_{i-1}^* D_{2s}^* \Gamma_{i-1} C_{2s} + \Lambda_i^* S_i + C_{2s}^* \Gamma_{i-1} C_{2s}\\
&=A_i^* K_i + K_i A_i + C_i^* K_i C_i + \Lambda_i^* \left[ R_i + D_{1s} \overline{P}_{11}^i D_{1s} \right]\Lambda_i\\
&\quad +\Lambda_i^* \left[ B_s^* \overline{P}_{11}^{i} + D^*_{1s} \overline{P}_{11}^i C_i + R_i \Theta_i + S_i \right]+ \left[ B_s^* \overline{P}_{11}^{i} + D^*_{1s} \overline{P}_{11}^i C_i + R_i \Theta_i + S_i \right]^* \Lambda_i\\
&\quad + \left[C_{2s} + D_{2s} \Theta_{i-1}\right]^* \Gamma_{i-1} \left[C_{2s} + D_{2s} \Theta_{i-1}\right].
\end{split}
\end{equation}
However, we observe that
\begin{equation}
\begin{split}
 B_s^* \overline{P}_{11}^{i} + D^*_{1s} \overline{P}_{11}^i C_i + R_i \Theta_i + S_i =  B_s^* \overline{P}_{11}^{i} + D^*_{1s} \overline{P}_{11}^i C_{1s} +S_i +\left[ R_i + D_{1s}^* \overline{P}_{11}^i D_{1s} \right]  \Theta_i = 0.
\end{split}
\end{equation}
Thus, we rewrite \eqref{temp4} as
\begin{equation}\label{temp5}
\begin{split}
-\frac{dK_i}{dt}  &=A_i^* K_i + K_i A_i + C_i^* K_i C_i + \Lambda_i^* \left[ R_s + D_{1s}^* \overline{P}_{11}^i D_{1s} + D_{2s}^* h_2(\overline{P}_{11}^{i-1}) D_{2s} \right]\Lambda_i\\
&\quad + \left[C_{2s} + D_{2s} \Theta_{i-1}\right]^* \Gamma_{i-1} \left[C_{2s} + D_{2s} \Theta_{i-1}\right]\\
&=A_i^* K_i + K_i A_i + C_i^* K_i C_i + \Lambda_i^* \left[ R_s + D_{1s}^* \overline{P}_{11}^{i+1} D_{1s} + D_{2s}^* h_2(\overline{P}_{11}^{i-1}) D_{2s} \right]\Lambda_i\\
&\quad + \Lambda_i^* D_{1s}^* K_i D_{1s} \Lambda_i + \left[C_{2s} + D_{2s} \Theta_{i-1}\right]^* \Gamma_{i-1} \left[C_{2s} + D_{2s} \Theta_{i-1}\right].
\end{split}
\end{equation}
To show that $K_i \geq 0$, for $i = 0,1,2,\ldots$, we will use mathematical induction. For $i = 0$, we have that
\begin{equation}
\begin{split}
-\frac{dK_0}{dt} &= A_s^* K_0 + K_0 A_s + C_{1s}^* K_0 C_{1s} + C_{2s}^*\left[ h_2'(\overline{P}_{11}^0) -  h_2(\overline{P}_{11}^0) \right] C_{2s}\\
&\quad+ \left[B_s^* \overline{P}_{11}^0 + D_{1s}^* \overline{P}_{11}^0 C_{1s} + D_{2s}^* h_2(\overline{P}_{11}^0) C_{2s} + L_s \right]^* \\
&\qquad\left(R_s + D_{1s}^* \overline{P}_{11}^0 D_{1s} +  D_{2s}^* h_2(\overline{P}_{11}^0)  D_{2s}\right)^{-1}\\
&\qquad \qquad \left[B_s^* \overline{P}_{11}^0 + D_{1s}^* \overline{P}_{11}^0 C_{1s} + D_{2s}^* h_2(\overline{P}_{11}^0) C_{2s} + L_s \right]\\
&\geq  A_s^* K_0 + K_0 A_s + C_{1s}^* K_0 C_{1s}.
\end{split}
\end{equation}
The above inequality comes from \eqref{lyapunovinq0.56} and the fact that $h_2'(P) \geq h_2(P)$ for all $P \geq 0$, see Lemma \ref{lemma: AREproperties}. Applying Theorem 2.1 of \cite{wonham1968matrix}, gives $K_0 \geq 0$. For $i = 1$, we can write \eqref{temp5} as
\begin{equation}
\begin{split}
\frac{dK_1}{dt}  + A_1^* K_1 + K_1 A_1 + C_1^* K_1 C_1 + \Lambda_1^* D_{1s}^* K_1 D_{1s} \Lambda_1\\
+ \Lambda_1^* \left[ R_s + D_{1s}^* \overline{P}_{11}^{2} D_{1s} + D_{2s}^* h_2(\overline{P}_{11}^{1}) D_{2s} \right]\Lambda_1 = 0.
\end{split}
\end{equation}
From \eqref{convexiteration}, we have that
\begin{equation}
-\frac{dK_1}{dt}  \geq A_1^* K_1 + K_1 A_1  + C_1^* K_1 C_1 + \Lambda_1^* D_{1s}^* K_1 D_{1s} \Lambda_1.
\end{equation}
By Theorem 2.1 of \cite{wonham1968matrix}, we get $K_1 \geq 0$. Now assume that $K_{n} \geq 0$ for $n = 0,1,2,\ldots, i-1$. Then $K_{i-2} = \overline{P}_{11}^{i-2} - \overline{P}_{11}^{i-1} \geq 0$ along with the fact that $h$ is increasing implies that $\Gamma_{i-1} \geq 0$. Thus, along with the inequality \eqref{convexiteration}, we have that
\begin{equation}
-\frac{dK_i}{dt}  \geq A_i^* K_i + K_i A_i + C_i^* K_1 C_i + \Lambda_i^* D_{1s}^* K_i D_{1s} \Lambda_i.
\end{equation}
Again by applying Theorem 2.1 of \cite{wonham1968matrix}, gives $K_i \geq 0$. Hence $\{\overline{P}_{11}^i\}_{i = 0}^\infty$ is a decreasing sequence
\begin{equation}
\overline{P}^0_{11}(t) \geq \overline{P}^i_{11}(t) \geq \overline{P}^{i+1}_{11}(t) \geq \alpha I,\quad \forall t\in [0,T],\ i = 0,1,2,\ldots
\end{equation}
where $\alpha = \inf_{k = 0,1,2,\ldots} \alpha_i \geq 0$. Thus the sequence $\{\overline{P}_{11}^i\}_{i = 0}^\infty$ is uniformly bounded. 

Denote $a = \sup_{t\in [0,T]} |K_0(t)|.$ From \eqref{temp5} and $K_i(T) = 0$, we can write, for $i = 1,2,3,\ldots$,
\begin{equation}
\begin{split}
K_i(t) &= \int_t^T \left[A_i^* K_i + K_i A_i + \widetilde{C}_i^* K_i \widetilde{C}_i + \Lambda_i^* \left[ R_s + D_{1s}^* \overline{P}_{11}^{i+1} D_{1s} + D_{2s}^* h_2(\overline{P}_{11}^{i-1}) D_{2s} \right]\Lambda_i \right.\\
&\left.\quad + \left[C_{2s} + D_{2s} \Theta_{i-1}\right]^* \Gamma_{i-1} \left[C_{2s} + D_{2s} \Theta_{i-1}\right] \right] ds.
\end{split}
\end{equation}
From now till the end of the proof, we shall denote by $M_T$ any positive constant, which may depend on $T$, uniformly in $i$. Observe that $|\Lambda_i|\leq M_T\left[ |K_{i-1}| + |\Gamma_{i-1}| \right]$, and since $h$ is Lipschitz, we also have that $|\Gamma_{i-1}|\leq M_T|K_{i-2}|$. As a result, we have that for $i = 1$
\begin{equation}
|K_1(t)|\leq M_T \int_t^T \left[ |K_1(s)| + |K_0(s)| \right] ds,
\end{equation}
and Gronwall's inequality implies that
\begin{equation}
|K_1(t)|\leq aM_T (T-t),
\end{equation}
For $i = 2,3,4,\ldots$, we have that
\begin{equation}
|K_i(t)|\leq M_T \int_t^T \left[ |K_i(s)| + |K_{i-1}(s)| + |K_{i-2}(s)| \right] ds.
\end{equation}
This implies that for $i = 2,3,4,\ldots$
\begin{equation}\label{finaleq}
|K_i(t)| + |K_{i-1}(t)|\leq M_T \int_t^T \left[ |K_i(s)| + |K_{i-1}(s)| + |K_{i-2}(s)| + |K_{i-3}(s)| \right] ds.
\end{equation}
Let $G_i(t) = |K_i(t)| + |K_{i-1}(t)|$. We will focus only on the odd indices $i = 1,3,5,\ldots$. We can calculate the upper bound of $G_1(t)$ explicitly as
\begin{equation}
G_1(t) = |K_1(t)| + |K_0(t)| \leq aM_T(T-t) + a.
\end{equation}
We can rewrite \eqref{finaleq} as
\begin{equation}
G_i(t)\leq M_T \int_t^T \left[ G_i(s) + G_{i-2}(s)\right] ds,\quad i = 3,5,7,\ldots,
\end{equation}
and applying Gronwall's inequality gives
\begin{equation}
G_i(t)\leq M_T \int_t^T  G_{i-2}(s) ds.
\end{equation}
By induction, we can see that for $k = 0,1,2,\ldots$,
\begin{equation}
G_{2k+1}(t) \leq a \left[ \frac{M_T^{k+1} (T-t)^{k+1}}{(k+1)!} +  \frac{M^{k} (T-t)^{k}}{k!}  \right].
\end{equation}
So $G_{2k+1}(t) = |K_{2k+1}(t)| + |K_{2k}(t)|$ implies that
\begin{equation}
|K_{2k+1}(t)|, |K_{2k}(t)| \leq a \left[ \frac{M_T^{k+1} (T-t)^{k+1}}{(k+1)!} +  \frac{M^{k} (T-t)^{k}}{k!}  \right].
\end{equation}
Hence the sequence $\{\overline{P}_{11}^i\}_{i = 1}^\infty$ is uniformly convergent, and we denote the limit as $\overline{P}_{11}^+$. Since $h_2$ is continuously differentiable, we have that for almost every $t\in [0,T]$
\begin{equation}
R_s + D_{1s}^* \overline{P}_{11}^+(t) D_{1s} +  D_{2s}^* h_2(\overline{P}_{11}^+(t))  D_{2s} = \lim_{i \rightarrow \infty}  R_s + D_{1s}^* \overline{P}_{11}^{i+1}(t) D_{1s} +  D_{2s}^* h_2(\overline{P}_{11}^{i-1}(t))  D_{2s} \geq \lambda I,
\end{equation}
and
\begin{equation}
\begin{split}
\begin{cases}
\lim_{i\rightarrow \infty}\Gamma_i &= 0,\\
\lim_{i\rightarrow \infty} \Theta_i &= - \left(R_s + D_{1s}^* \overline{P}_{11}^+ D_{1s} +  D_{2s}^* h_2(\overline{P}_{11}^+)  D_{2s}\right)^{-1}\\
&\qquad \qquad \left[B_s^* \overline{P}_{11} + D_{1s}^* \overline{P}_{11} C_{1s} + D_{2s}^* h_2(\overline{P}_{11}^+) C_{2s} + L_s \right],\\
\lim_{i\rightarrow \infty}A_i &= A_s + B_s \Theta,\quad \quad \lim_{i\rightarrow \infty}C_i = C_{1s} + D_{1s} \Theta,\\
\lim_{i\rightarrow \infty}Q_i &= Q_s + C_{2s}^* h_2(\overline{P}_{11}^+) C_{2s},\quad \lim_{i\rightarrow \infty}R_i = R_s + D_{2s}^* h_2(\overline{P}_{11}^+)  D_{2s},\\
\lim_{i\rightarrow \infty}S_i &= L_s + D_{2s}^* h_2(\overline{P}_{11}^+) C_{2s}.
\end{cases}
\end{split}
\end{equation}
Moreover, $\overline{P}_{11}^+ = \lim_{i\rightarrow \infty} \overline{P}_{11}^i \geq \alpha I$. Hence $\overline{P}_{11}^+$ is the solution of the Lyapunov equation
\begin{equation}\label{finallyapunov}
\begin{split}
\begin{cases}
\frac{d\overline{P}_{11}}{dt} + \left[A_s + B_s \Theta\right]^* \overline{P}_{11} + \overline{P}_{11} \left[A_s + B_s \Theta \right] + \left[C_{1s} + D_{1s}\Theta \right]^* \overline{P}_{11} \left[C_{1s} + D_{1s}\Theta \right] \\
\ \qquad + \Theta^* \left[R_s + D_{2s}^* h_2(\overline{P}_{11})  D_{2s}\right] \Theta + \left[L_s + D_{2s}^* h_2(\overline{P}_{11}) C_{2s}\right]^* \Theta + \Theta^*\left[L_s + D_{2s}^* h_2(\overline{P}_{11}) C_{2s}\right] \\
\ \qquad+ Q_s + C_{2s}^* h_2(\overline{P}_{11}) C_{2s} = 0,\\
\overline{P}_{11}(T) = 0.
\end{cases}
\end{split}
\end{equation}
It is then straightforward to show that \eqref{finallyapunov} is equivalent to the Riccati equation \eqref{reducedDREfixed}. Finally, setting $\overline{P}_{22}^+ = h_2(\overline{P}_{11}^+)$ and $\overline{F}_2^+ = \overline{F}_2(\overline{P}_{11}^+)$ completes the proof.
\qed 

\begin{corollary}\label{corollary: reducedsystem}
Suppose that Assumption \ref{mainassumptions} holds. Let $(\overline{P}_{11}^+,\overline{P}_{22}^+)$ be the solution to the reduced differential-algebraic Riccati equation \eqref{reducedDRE}-\eqref{reducedARE} and $\overline{F}_2^+$ the feedback operator shown in Theorem \ref{theorem: DARE}. Then the reduced system \eqref{reducedsystem} admits solution $(\overline{P}_{11}, \overline{P}_{12},\overline{P}_{22}) \in C([0,T]; \mathbb{S}^{n_1}_+) \times C([0,T]; \R^{n_1\times n_2}) \times C([0,T]; \mathbb{S}^{n_2}_{++})$ where $(\overline{P}_{11}, \overline{P}_{22}) = (\overline{P}_{11}^+, \overline{P}_{22}^+)$ and
\begin{equation}\label{P12soln}
\begin{split}
\overline{P}_{12} &=- \left[ \overline{P}_{11}^+ A_{12} + A_{21}^* \overline{P}_{22}^+ + C_{11}^* \overline{P}_{11}^+ C_{12} + C_{21}^* \overline{P}_{22}^+ C_{22} \right.\\
&\left. \qquad + Q_{12} + \left(  \overline{P}_{11}^+B_1  + C_{11}^* \overline{P}_{11}^+ D_1 + C_{21}^* \overline{P}_{22}^+ D_2 \right) \overline{F}_2^+ \right] 
\left[ A_{22} + B_2 \overline{F}_2^+ \right]^{-1}.
\end{split}
\end{equation} In addition, $\overline{P}_{22}$ is unique in the set of positive-semidefinite solutions. 
\end{corollary}
\noindent \textit{Proof.} The result follows from Theorem \ref{theorem: equiv}, Lemma \ref{lemma: AREproperties} and Theorem \ref{theorem: DARE}.
\qed

\section{Convergence and the Tikhonov Theorem}\label{section: tikhonov}

Let $p > 0$ be an arbitrary constant such that 
\begin{equation}
|\overline{P}_{11}^+(t)|  < p,\quad \forall t\in [0,T].
\end{equation}
Fix $(t,P_{11}) \in [0,T] \times B_p$ where $B_p = \{P_{11} \in \mathbb{S}^{n_1}_+\ |\ |P_{11}| \leq p \}$. We work with the "stretched" time variable $\tau = t/\epsilon$. Consider the boundary-layer problem
\begin{subnumcases}{\label{boundary}}
\label{boundary1}
\frac{d\widehat{P}_{12}}{d\tau} = g_1(P_{11}, \widehat{P}_{12} + h_1(P_{11}), \widehat{P}_{22} + h_2(P_{11}),0),\\
\label{boundary2}
\frac{d\widehat{P}_{22}}{d\tau} = g_2(P_{11}, \widehat{P}_{12} + h_1(P_{11}), \widehat{P}_{22} + h_2(P_{11}),0),
\end{subnumcases}
with initial conditions $(\overline{P}_{12}(0), \overline{P}_{22}(0))$ where $(h_1(P_{11}), h_2(P_{11}))$ is the solution to the reduced system \eqref{reduced2}-\eqref{reduced3} when $P_{11}$ is fixed. More precisely, $h_2(P_{11})$ is as defined in Lemma \ref{lemma: AREproperties} with the property that the eigenvalues of $A_{22} + B_2 \overline{F}_2$ have negative real parts where
\begin{equation}\label{temp9}
\begin{split}
\overline{F}_2 &= -\left[R + D_1^* P_{11} D_1 + D_2^*  h_2(P_{11}) D_2\right]^{-1} \left[B_2^*  h_2(P_{11}) + D_1^* P_{11} C_{12} + D_2^*  h_2(P_{11}) C_{22}\right].
\end{split}
\end{equation}
In addition, $h_1$ is defined as the solution $\overline{P}_{12}$ to the reduced system in Corollary \ref{corollary: reducedsystem} but parametrised by $P_{11}$. That is,
\begin{equation}
\begin{split}
h_1(P_{11}) &= -\left[ P_{11} A_{12} + A_{21}^* h_2(P_{11})  + C_{11}^* P_{11} C_{12} + C_{21}^* h_2(P_{11}) C_{22} + Q_{12} \right.\\
&\left.\quad+ \left( P_{11} B_1 + C_{11}^* P_{11} D_1+ C_{21}^*  h_2(P_{11}) D_2 \right) \overline{F}_2 \right] \left[ A_{22} + B_2 \overline{F}_2 \right]^{-1}.
\end{split}
\end{equation}
When $\epsilon$ tends towards zero, it is equivalent to the variable $\tau$ tending towards infinity. For this reason, we are interested in the equilibrium $(0,0)$ of the boundary-layer problem \eqref{boundary}.

For $P_{12} \in \R^{n_1\times n_2}$ and $P_{22} \in \mathbb{S}^{n_2}$, let us define the norm of the pair $(P_{12}, P_{22})$ as $|(P_{12}, P_{22})| := |P_{12}| + |P_{22}|$. Let $q_0$ be a positive constant, independent of $p$, and define the set $B_{q_0}$ of initial values as
\begin{equation}
\begin{split}
B_{q_0} := \{ (\widehat{P}_{12}(0), \widehat{P}_{22}(0)) \in \R^{n_1\times n_2} \times \mathbb{S}^{n_2}\ |\ h_2(0) + \widehat{P}_{22}(0) \geq 0 \text{ and } |(\widehat{P}_{12}(0), \widehat{P}_{22}(0))| < q_0 \}.
\end{split}
\end{equation}
In Theorem \ref{theorem: tikhonov}, we will be using the version of the boundary-layer problem \eqref{boundary} when $P_{11} = 0$ with initial values $(\widehat{P}_{12}(0),\widehat{P}_{22}(0)) = (-h_1(0),-h_2(0))$. For this reason, we will be forward looking and set $q_0$ large enough such that $|(-h_1(0),-h_2(0))| < q_0$. This implies that $(-h_1(0),-h_2(0)) \in B_{q_0}$. 

\begin{lemma}\label{lemma: boundarylayerproblem}
Suppose that Assumption \ref{mainassumptions} holds. Then the boundary-layer problem \eqref{boundary} with initial values $(\widehat{P}_{12}(0),\widehat{P}_{22}(0)) \in B_{q_0}$ has a solution $(\widehat{P}_{12}, \widehat{P}_{22}) \in C([0,T]; \R^{n_1\times n_2})\times C([0,T]; \mathbb{S}^{n_2})$ satisfying $h_2(P_{11}) + \widehat{P}_{22}(\tau)\geq 0,\forall \tau \geq 0$ and converges to the equilibrium $(0,0)$ as $\tau \rightarrow \infty$.  Moreover, the equilibrium $(0,0)$ is exponentially stable, uniformly in $P_{11}$. That is, for all $(\widehat{P}_{12}(0),\widehat{P}_{22}(0)) \in B_{q_0}$,
\begin{equation}\label{expstab}
\begin{split}
\begin{cases}
|\widehat{P}_{12}(\tau)| \leq \widehat{M} e^{-\gamma \tau} \left[|\widehat{P}_{12}(0)| + |\widehat{P}_{22}(0)|\right],\quad \forall P_{11} \in B_p,\ \forall \tau \geq 0,\\
|\widehat{P}_{22}(\tau)| \leq \widehat{M} e^{-\frac{3\gamma}{2} \tau} |\widehat{P}_{22}(0)|,\quad \forall P_{11} \in B_p,\ \forall \tau \geq 0,
\end{cases}
\end{split}
\end{equation}
for some positive constants $\widehat{M}$, which may depend on $p$ and $q_0$, and $\gamma > 0$ is as defined in \eqref{expstabG}.
\end{lemma}
\noindent \textit{Proof.} Let $Z_1(\tau) = \widehat{P}_{12}(\tau) + h_1(P_{11})$ and $Z_2(\tau) = \widehat{P}_{22}(\tau) + h_2(P_{11})$ for $\tau \geq 0$. The boundary-layer problem \eqref{boundary2} can be restated as
\begin{equation}\label{boundary2temp}
\frac{dZ_2}{d\tau} = g_2(P_{11}, Z_1, Z_2,0),\quad Z_2(0) = h_2(P_{11}) + \widehat{P}_{22}(0).
\end{equation}
It is clear that \eqref{boundary2temp} is equivalent to the Riccati equation \eqref{reducedAREfixeddiff} with initial value $h_2(P_{11})+\widehat{P}_{22}(0)$. Observe that, from Lemma \ref{lemma: AREproperties}, $h_2(P_{11})+\widehat{P}_{22}(0) \geq 0$. Moreover, the equation \eqref{boundary2temp} admits a solution $Z_2 \in C([0,T]; \mathbb{S}^{n_2}_+)$ that converges to $h_2(P_{11}) > 0$, uniformly in $P_{11} \in B_p$ and moreover,
\begin{equation}\label{temp10a}
|Z_2(\tau) - h_2(P_{11})| \leq M_P e^{-\frac{3\gamma}{2} \tau} |\widehat{P}_{22}(0)|,\quad \forall P_{11} \in B_p,\ \forall \tau \geq 0,
\end{equation}
where $M_P$ and $\gamma$ are positive constants defined in \eqref{expstabG}. We note that $M_P$ may depend on $p$. Hence, we have that
\begin{equation}\label{temp10}
|\widehat{P}_{22}(\tau)| \leq M_P e^{-\frac{3\gamma}{2} \tau} |\widehat{P}_{22}(0)|,\quad \forall P_{11} \in B_p,\ \forall \tau \geq 0.
\end{equation}
In addition, Lemma \ref{lemma: AREproperties} tells us that the eigenvalues of $A_{22} + B_2 \overline{F}_2$ have negative real parts where
\begin{equation}\label{temp8}
\begin{split}
\overline{F}_2(P_{11}) &= - (R + D_1^* P_{11} D_1 + D_2^* h_2(P_{11}) D_2)^{-1} \\
&\qquad \qquad\left[B_2^* h_2(P_{11}) + D_1^* P_{11} C_{12} + D_2^* h_2(P_{11}) C_{22}\right].
\end{split}
\end{equation}
and for some positive constant $M$,
\begin{equation}\label{expstab100}
e^{\tau(A_{22} + B_2 \overline{F}_2)} \leq M e^{-\gamma\tau},\quad \tau \geq 0.
\end{equation}
We point out that the constant $\gamma$ in \eqref{temp10} and \eqref{expstab100} are the same.

Next, consider the equation \eqref{boundary1}. From the formulation \eqref{reduced5}, we can show that
\begin{equation}
g_1(P_{11}, \widehat{P}_{12} + h_1(P_{11}), \widehat{P}_{22} + h_2(P_{11}),0) = \widehat{P}_{12}\left[ A_{22} + B_2 \widetilde{F}_2(P_{11})\right] + L(P_{11},\widehat{P}_{22})
\end{equation}
where $L(P_{11},\widehat{P}_{22})$ is defined as
\begin{equation}
\begin{split}
&L(P_{11},\widehat{P}_{22}) \\
&\quad= A_{21}^* \widehat{P}_{22} + C_{21}^* \widehat{P}_{22} C_{22} + C_{21}^* \widehat{P}_{22} D_2 \widetilde{F}_2(P_{11}) \\
&\qquad + \left[ P_{11} B_1 + h_1(P_{11}) B_2 + C_{11}^* P_{11} D_1 + C_{21}^* h_2(P_{11}) D_2 \right] \left[ \widetilde{F}_2(P_{11}) - \overline{F}_2(P_{11})\right]
\end{split}
\end{equation}
and $\widetilde{F}_2(P_{11})$ is defined as
\begin{equation}
\begin{split}
\widetilde{F}_2(P_{11}) &= -\left[R + D_1^* P_{11} D_1 + D_2^*  \left(\widehat{P}_{22} + h_2(P_{11}) \right) D_2\right]^{-1}\\
&\quad \left[B_2^*  \left(\widehat{P}_{22} + h_2(P_{11}) \right) + D_1^* P_{11} C_{12} + D_2^*  \left(\widehat{P}_{22} + h_2(P_{11}) \right) C_{22}\right].
\end{split}
\end{equation}
Note that $R + D_1^* P_{11} D_1 + D_2^*  \left(\widehat{P}_{22} + h_2(P_{11}) \right) D_2$ is indeed invertible as $R$ is strictly positive, and $P_{11}$ and $\widehat{P}_{22}(\tau) + h_2(P_{11})$ are both positive semidefinite for all $\tau \geq 0$.

Now consider the matrix $A_{22} + B_2 \widetilde{F}_2(P_{11}) = A_{22} + B_2 \overline{F}_2 + B_2 \left[ \widetilde{F}_2(P_{11})- \overline{F}_2(P_{11})\right]$. Since the eigenvalues of $G := A_{22} + B_2 \overline{F}_2$ have negative real parts, it is generator of an exponentially stable $C_0$-semigroup $S_G$ satisfying the inequality \eqref{expstab100}. Fix $\tau \geq 0$. Let $U(\tau,s), 0\leq s\leq \tau$ be the solution to the backwards equation
\begin{equation}
\begin{split}
\begin{cases}
\frac{d}{ds} U(\tau,s) = -(A_{22} + B_2 \widetilde{F}_2(P_{11})(s)) U(\tau,s),\\
U(\tau,\tau) = I,
\end{cases}
\end{split}
\end{equation}
in which $U(\tau,s)$ satisfies
\begin{equation}\label{evolutionfam}
U(\tau,s) = e^{(\tau-s)(A_{22} + B_2 \overline{F}_2)} + \int_s^\tau e^{(r - s)(A_{22} + B_2 \overline{F}_2)} B_2\left[ \widetilde{F}_2(P_{11})(r)- \overline{F}_2(P_{11})\right] U(\tau,r) dr.
\end{equation}
This comes from the fact that $(\overline{U}_G(\tau,s)^*)_{0\leq s\leq \tau}$ is an evolution operator with the generator $G^*$, see \cite{curtain1976infinite}. Moreover, from \eqref{expstab100} and \eqref{evolutionfam}, we have that for any $0\leq s\leq \tau$,
\begin{align*}
|U(\tau,s)| \leq M e^{-\gamma(\tau-s)} + M |B_2| \int_s^\tau e^{-\gamma(\tau-r)} | \widetilde{F}_2(P_{11})(r)- \overline{F}_2(P_{11})|\ |U(t,r)|dr.
\end{align*}
From \eqref{temp9} and \eqref{temp8}, we notice that
\begin{equation}\label{lipschitzzz}
|\widetilde{F}_2(P_{11})(r) - \overline{F}_2(P_{11})| \leq M_1 |\widehat{P}_{22}(r)|,
\end{equation}
where $M_1$ is a positive constant that depends on $p$. Thus, from \eqref{temp10}, we have that
\begin{align*}
|U(\tau,s)| &\leq M e^{-\gamma(\tau-s)} + M M_1 |B_2| \int_s^\tau e^{-\gamma(r - s)} |\widehat{P}_{22}(r)| |U(\tau,r)|dr\\
&\leq M e^{-\gamma(\tau-s)} + M M_1 M_G |B_2|q_0 \int_s^\tau e^{-\gamma(r - s)} e^{-\frac{3\gamma}{2} r} |U(\tau,r)|dr.
\end{align*}
Hence
\begin{align*}
e^{- \gamma s} |U(\tau,s)| 
&\leq M e^{-\gamma \tau} + M M_1 M_G |B_2| q_0 \int_s^\tau e^{-\frac{3\gamma}{2} r} \left[e^{-\gamma r} |U(t,r)|\right]dr.
\end{align*}
By Gronwall's inequality, we have that $U(\tau,s), 0 \leq s\leq \tau$ satisfies the inequality
\begin{equation}\label{tempexp}
|U(\tau,s)| \leq M e^{-\gamma (\tau - s)} e^{\frac{2M M_1 M_G |B_2|q_0}{3\gamma}}.
\end{equation}
Now, for fixed $\tau \geq 0$, let us differentiate the mapping
\begin{equation}
[0,\tau] \ni s \mapsto \widehat{P}_{12}(s) U(\tau,s),
\end{equation}
to obtain
\begin{align*}
\frac{d}{ds}\left[\widehat{P}_{12}(s) U(\tau,s)\right] &= \widehat{P}_{12}(s)\left[ A_{22} + B_2 \widetilde{F}_2(P_{11})(s)\right] U(\tau,s) + L(P_{11},\widehat{P}_{22}(s)) U(\tau,s) \\
&\quad - \widehat{P}_{12}\left[ A_{22} + B_2 \widetilde{F}_2(P_{11})(s)\right] U(\tau,s)\\
&= L(P_{11},\widehat{P}_{22}(s)) U(\tau,s).
\end{align*}
Integrating from $0$ to $\tau$ gives
\begin{align*}
\widehat{P}_{12}(\tau) = \widehat{P}_{12}(0)U(\tau,0) + \int_0^\tau L(P_{11},\widehat{P}_{22}(s))U(\tau,s) ds.
\end{align*}
From the definition of $L$, \eqref{lipschitzzz} and \eqref{temp10} there exists a positive constant $M_L$, depending on $p$ such that
\begin{equation}
| L(P_{11}, \widehat{P}_{22}(s)) | \leq M_L |\widehat{P}_{22}(s)| \leq M_L M_P e^{-\frac{3\gamma}{2} s} |\widehat{P}_{22}(0)|, \quad s\geq 0.
\end{equation} 
Hence, along with \eqref{tempexp}, we have that
\begin{align*}
|\widehat{P}_{12}(\tau)| &\leq M e^{\frac{2M M_1 M_G |B_2|q_0}{3\gamma}} e^{-\gamma \tau}  |\widehat{P}_{12}(0)| + M M_L M_P e^{\frac{2M M_1 M_G |B_2|q_0}{3\gamma}} \int_0^\tau e^{-\frac{3\gamma}{2} s} e^{-\gamma (\tau-s)} |\widehat{P}_{22}(0)| ds\\
&\leq M e^{\frac{2M M_1 M_G |B_2|q_0}{3\gamma}} e^{-\gamma \tau}  |\widehat{P}_{12}(0)| + \frac{2M M_L M_P}{\gamma} e^{\frac{2M M_1 M_G |B_2|q_0}{3\gamma}} e^{-\gamma \tau} |\widehat{P}_{22}(0)|.
\end{align*}
\qed

Next, let $q$ be an arbitrarily large positive constant such that $q > (3 \widehat{M} + 1) q_0$ and define the set $B_q$ as
\begin{equation}
\begin{split}
B_{q} := \{ (\widehat{P}_{12}, \widehat{P}_{22}) \in \R^{n_1\times n_2} \times \mathbb{S}^{n_2}\ |\ h_2(P_{11}) + \widehat{P}_{22} \geq 0, \forall P_{11} \in B_p  \text{ and } |(\widehat{P}_{12}, \widehat{P}_{22})| < q\}.
\end{split}
\end{equation}
Note that $q$ may depend on both $p$ and $q_0$. Here we make a couple of observations. First, from Lemma \ref{lemma: boundarylayerproblem}, the solution to the boundary-layer problem \eqref{boundary} is contained in $B_q$. Second, since $h_2$ is an increasing function from Lemma \ref{lemma: AREproperties}, we have that $B_{q_0} \subset B_q$.

\begin{theorem}\label{theorem: tikhonov}
Suppose that Assumption \ref{mainassumptions} hold. Let $T > 0$ be any finite time horizon and let
\begin{itemize}
\item $(P^\epsilon_{11}, P^\epsilon_{12}, P^\epsilon_{22})$ be the unique solution to the full system \eqref{fullsystem};
\item $(\overline{P}_{11}, \overline{P}_{12}, \overline{P}_{22})$ be the solution to the reduced system \eqref{reducedsystem} defined in Corollary \ref{corollary: reducedsystem};
\item $(\widehat{P}_{12},\widehat{P}_{22})$ be the solution to the boundary-layer problem \eqref{boundary}.
\end{itemize}
Then there exists a positive constant $\epsilon^*$ such that for all $0 < \epsilon< \epsilon^*$
\begin{equation}
\begin{split}
\begin{cases}
P^{\epsilon}_{11}(t) - \overline{P}_{11}(t) = O(\epsilon),\\
P^{\epsilon}_{12}(t) - \overline{P}_{12}(t) - \widehat{P}_{12}\left(\frac{T-t}{\epsilon}\right) = O(\epsilon),\\
P^{\epsilon}_{22}(t) - \overline{P}_{22}(t) - \widehat{P}_{22}\left(\frac{T-t}{\epsilon}\right) = O(\epsilon),
\end{cases}
\end{split}
\end{equation}
uniformly in $t\in [0,T]$.
\end{theorem}

\noindent \textit{Proof.} The result follows from Corollary \ref{corollary: reducedsystem}, Lemma \ref{lemma: boundarylayerproblem} and the Tikhonov Theorem, see Theorem 9.1 of \cite{khalil1996nonlinear}. In order for the Tikhonov Theorem to apply, we need to make two observations.

First, referencing the notation used in Theorem 9.1 of \cite{khalil1996nonlinear}, we mention that in order for the Tikhonov Theorem to hold, we need make sure that the norm of initial values $(-h_1(0),-h_2(0))$ of the boundary-layer problem \eqref{boundary} are bounded by a particular constant $\mu > 0$, which is directly proportional to the constant $q_0 > 0$. This is resolved by observing that we can make $q_0$ arbitrarily large, which ensures that $|(-h_1(0),-h_2(0))| < \mu$.

Second, note that we have included a positivity condition in the sets $B_p$ and $B_q$, where its product $B_p \times B_q$ can be described as the intersection between the region
\begin{equation}
\mathcal{R} = \{ (P_{11}, P_{12}, P_{22}) \in \mathbb{S}^{n_1}\times \R^{n_1\times n_2} \times \mathbb{S}^{n_2}\ |\ P_{11} \geq 0 \text{ and } h_2(P_{11}) + P_{22} \geq 0\},
\end{equation}
and open and closed balls. From Lemma \ref{lemma: fullsystemsoln} and Corollary \ref{corollary: reducedsystem} the solutions of the full system \eqref{fullsystem} and reduced system \eqref{reducedsystem} both belong to the set $\mathcal{R}$ for all $t\in [0,T]$. From Lemma \ref{lemma: boundarylayerproblem}, the same can be said about the solution of the boundary-layer problem \eqref{boundary} as $q$ can be chosen arbitrarily large. Hence all processes that we consider are contained in the region $\mathcal{R}$ and thus, we can work with the sets $B_p$ and $B_{q}$ and treat them as balls. The same can be said for similar sets which are constructed by an intersection of $\mathcal{R}$ and open or closed balls.
\qed

The follow result describes the convergence of the solution to the Riccati equation \eqref{riccati}.

\begin{corollary}
Suppose that Assumption \ref{mainassumptions} holds. Let $T > 0$ be any finite time horizon and $\epsilon^*$ be the small positive parameter defined in Theorem \ref{theorem: tikhonov}. Let $P^\epsilon \in C([0,T]; \mathbb{S}^{n}_+)$ be the solution to the Riccati equation \eqref{riccati}. Then for all $0 < \epsilon< \epsilon^*$, we have that
\begin{equation}
P^\epsilon(t) - \begin{pmatrix}
\overline{P}_{11} & 0\\
0 & 0
\end{pmatrix} = O(\epsilon)
\end{equation}
uniformly in $t\in [0,T]$.
\end{corollary}

\noindent \textit{Proof}. The result follows from the first-order representation \eqref{partition}, Theorem \ref{theorem: tikhonov} and the uniform boundedness of the solutions to the reduced system and boundary-layer problem with respect to time.
\qed

\begin{corollary}\label{corollary: integral}
Suppose that Assumption \ref{mainassumptions} holds. Let $T > 0$ be any finite time horizon and $\epsilon^*$ be the small positive parameter defined in Theorem \ref{theorem: tikhonov}. Then for any positive integer $j$, there exists a positive constant $K(T,j)$ such that for all $0 < \epsilon < \epsilon^*$ and $i = 1,2$
\begin{equation}
\int_0^T |P^\epsilon_{i2}(t) - \overline{P}_{i2}(t)|^j dt \leq \epsilon K(T,j).
\end{equation}
\end{corollary}

\noindent \textit{Proof.} Fix $T > 0$ and let $t \in [0,T]$. For $i = 1,2$, Lemma \ref{lemma: boundarylayerproblem} and Theorem \ref{theorem: tikhonov} imply that
\begin{align*}
|P^\epsilon_{i2}(t) - \overline{P}_{i2}(t)| &\leq \epsilon K_1(T) + \left| \widehat{P}_{i2}\left(\frac{T-t}{\epsilon}\right)\right|\\
&\leq \epsilon K_1(T) + K_2 e^{-\frac{\gamma (T-t)}{\epsilon}}
\end{align*}
for some positive constants $K_1(T)$, which depends on $T$, and $K_2$. Hence for all positive integers $j$, we have
\begin{align*}
|P^\epsilon_{i2}(t) - \overline{P}_{i2}(t)|^j &\leq  2^{j-1}\left(\epsilon^j K_1(T)^j + K_2^j e^{-\frac{j\gamma (T-t)}{\epsilon}}\right).
\end{align*}
Finally, integrating the above inequality gives the desired result.
\qed
\section{Approximately optimal control and estimation of the value function}\label{section: estimates}

In this section, we use Theorem \ref{theorem: tikhonov} to construct an approximate optimal control and value function based on the solution to the reduced system \eqref{reducedsystem}. We preface by stating that the letter $K$ will be reserved for a positive constant and is not necessarily the same in each instance. In the situations where $K$ may depend on another relevant constant, say $T$, we will denote this as $K(T)$.

Recall from Theorem \ref{theorem: optimality}, the optimal control $\widehat{u}^\epsilon$ is given by
\begin{equation}
\widehat{u}^\epsilon(t) = \widehat{F}^\epsilon(t) \widehat{X}^\epsilon(t) = \widehat{F}^\epsilon_1(t) \widehat{X}^\epsilon_1(t) + \widehat{F}^\epsilon_2(t) \widehat{X}^\epsilon_2(t),\quad \forall t\in [0,T],
\end{equation}
where the feedback operators are given by
\begin{subnumcases}{}
\begin{split}
\widehat{F}_1^\epsilon &= - (\Delta^\epsilon)^{-1} \left[ B_1^* P_{11}^\epsilon + B_2^* (P_{12}^\epsilon)^* + D_1^* P_{11}^\epsilon C_{11} \right.\\
&\left. \quad + \sqrt{\epsilon}\left( D_2^* (P_{12}^\epsilon)^* C_{11} + D_1^* P_{12}^\epsilon C_{21} \right) + D_2^* P_{22}^\epsilon C_{21}\right],
\end{split}\\
\begin{split}
\widehat{F}_2^\epsilon  &= - (\Delta^\epsilon)^{-1}  \left[  \epsilon B_1^* P_{12}^\epsilon + B_2^* P_{22}^\epsilon + D_1^* P_{11}^\epsilon C_{12} \right.\\
&\left. \quad + \sqrt{\epsilon}\left( D_2^* (P_{12}^\epsilon)^* C_{12} + D_1^* P_{12}^\epsilon C_{22} \right) + D_2^* P_{22}^\epsilon C_{22} \right],
\end{split}
\end{subnumcases}
and $\widehat{X}^\epsilon_1(t) := X_1(t; \widehat{u}^\epsilon)$ and $\widehat{X}^\epsilon_2(t) := X_2^\epsilon(t; \widehat{u}^\epsilon)$ are solutions to the optimal state equations
\begin{equation}
\begin{split}
\begin{cases}
d\widehat{X}^\epsilon_1(t) = \left[ \left( A_{11} + B_1 \widehat{F}_1^\epsilon(t) \right) \widehat{X}^\epsilon_1(t) + \left( A_{12} + B_1 \widehat{F}_2^\epsilon(t) \right) \widehat{X}^\epsilon_2(t)  \right]dt\\
\qquad \qquad + \left[ \left( C_{11} + D_1 \widehat{F}_1^\epsilon(t) \right) \widehat{X}^\epsilon_1(t) + \left( C_{12} + D_1 \widehat{F}_2^\epsilon(t) \right) \widehat{X}^\epsilon_2(t)  \right] dW(t),\\
\widehat{X}^\epsilon_1(0) = x_1,
\end{cases}
\end{split}
\end{equation}
and
\begin{equation}
\begin{split}
\begin{cases}
d\widehat{X}^\epsilon_2(t) = \frac{1}{\epsilon}  [\left[ \left( A_{21} + B_2 \widehat{F}_1^\epsilon(t) \right) \widehat{X}^\epsilon_1(t) + \left( A_{22} + B_2 \widehat{F}_2^\epsilon(t) \right) \widehat{X}^\epsilon_2(t)  \right]dt\\
\qquad \qquad + \frac{1}{\sqrt{\epsilon}} \left[ \left( C_{21} + D_2 \widehat{F}_1^\epsilon(t) \right) \widehat{X}^\epsilon_1(t) + \left( C_{22} + D_2 \widehat{F}_2^\epsilon(t) \right) \widehat{X}^\epsilon_2(t)  \right] dW(t),\\
\widehat{X}^\epsilon_2(0) = x_2.
\end{cases}
\end{split}
\end{equation}

We construct an approximately optimal control $\overline{u}^\epsilon$ by formally setting $\epsilon = 0$ and using the solution to the reduced system in the feedback operators. In doing so, we obtain
\begin{equation}\label{approxoptimcontrol}
\overline{u}^\epsilon(t) = \overline{F}_1(t) \overline{X}_1^\epsilon(t) + \overline{F}_2(t) \overline{X}_2^\epsilon(t)
\end{equation}
where
\begin{subnumcases}{}
\overline{F}_1 = - \overline{\Delta}^{-1} \left[ B_1^* \overline{P}_{11} + B_2^* \overline{P}_{12}^* + D_1^* \overline{P}_{11} C_{11} + D_2^* \overline{P}_{22} C_{21} \right],\\
\overline{F}_2  = - \overline{\Delta}^{-1}  \left[ B_2^* \overline{P}_{22} + D_1^* \overline{P}_{11} C_{12} + D_2^* \overline{P}_{22} C_{22}\right],
\end{subnumcases}
and $\overline{X}^\epsilon_1(t) := X_1(t; \overline{u}^\epsilon)$ and $\overline{X}^\epsilon_2(t) := X_2^\epsilon(t; \overline{u}^\epsilon)$ are solutions to the state equations
\begin{equation}\label{approxstate1}
\begin{split}
\begin{cases}
d\overline{X}^\epsilon_1(t) = \left[ \left( A_{11} + B_1 \overline{F}_1(t) \right) \overline{X}^\epsilon_1(t) + \left( A_{12} + B_1 \overline{F}_2(t) \right) \overline{X}^\epsilon_2(t)  \right]dt\\
\qquad \qquad + \left[ \left( C_{11} + D_1 \overline{F}_1(t) \right) \overline{X}^\epsilon_1(t) + \left( C_{12} + D_1 \overline{F}_2(t) \right) \overline{X}^\epsilon_2(t)  \right] dW(t),\\
\overline{X}^\epsilon_1(0) = x_1,
\end{cases}
\end{split}
\end{equation}
and
\begin{equation}\label{approxstate2}
\begin{split}
\begin{cases}
d\overline{X}^\epsilon_2(t) = \frac{1}{\epsilon}   \left[ \left( A_{21} + B_2 \overline{F}_1(t) \right) \overline{X}^\epsilon_1(t) + \left( A_{22} + B_2 \overline{F}_2(t) \right) \overline{X}^\epsilon_2(t)  \right]dt\\
\qquad \qquad + \frac{1}{\sqrt{\epsilon}} \left[ \left( C_{21} + D_2 \overline{F}_1(t) \right) \overline{X}^\epsilon_1(t) + \left( C_{22} + D_2 \overline{F}_2(t) \right) \overline{X}^\epsilon_2(t)  \right] dW(t),\\
\overline{X}^\epsilon_2(0) = x_2.
\end{cases}
\end{split}
\end{equation}
\begin{lemma}\label{lemma: feedbackbound}
Suppose that Assumption \ref{mainassumptions} holds. Let $T > 0$ be any finite time horizon and $\epsilon^*$ be the small positive parameter defined in Theorem \ref{theorem: tikhonov}. Then there exists a positive constant $K(T)$, which depends on $T$, such that for all $i = 1,2$ and $0 < \epsilon < \epsilon^*$
\begin{equation}
\int_0^T |\widehat{F}^\epsilon_i(t) - \overline{F}_i(t)|^2 dt \leq \epsilon K(T).
\end{equation}
\end{lemma}
\noindent \textit{Proof.} The result follows from Theorem \ref{theorem: tikhonov} and Corollary \ref{corollary: integral}. \qed

\begin{lemma}\label{lemma: statebound}
Suppose that Assumption \ref{mainassumptions} holds. Then for any finite time horizon $T > 0$, we have 
\begin{equation}
\sup_{\epsilon \in (0,1]} \sup_{t\in [0,T]} \E\left[ |\overline{X}_1^\epsilon(t)|^2 + |\overline{X}_2^\epsilon(t)|^2 \right] < \infty.
\end{equation}
\end{lemma}
\noindent \textit{Proof.} Fix $\epsilon \in (0,1]$. Let us begin with \eqref{approxstate1}. Taking the norm, we have that for all $0 \leq t\leq T$
\begin{align*}
|\overline{X}^\epsilon_1(t)| &\leq |x_1| + \int_0^t \left[ | A_{11} + B_1 \overline{F}_1(s) |\ |\overline{X}^\epsilon_1(s)| + | A_{12} + B_1 \overline{F}_2(s) | \ |\overline{X}^\epsilon_2(s)|  \right]ds\\
&\quad + \Bigg| \int_0^t \left[ \left( C_{11} + D_1 \overline{F}_1(s) \right) \overline{X}^\epsilon_1(s) + \left( C_{12} + D_1 \overline{F}_2(s) \right) \overline{X}^\epsilon_2(s)  \right] dW(s)\Bigg|.
\end{align*}
By Corollary \ref{corollary: reducedsystem}, the solution $(\overline{P}_{11},\overline{P}_{12},\overline{P}_{22})$ is bounded on $[0,T]$ and consequently $\overline{F}_1$ and $\overline{F}_2$ is as well. So we can write the above as
\begin{align*}
|\overline{X}^\epsilon_1(t)| &\leq |x_1| + K(T) \int_0^t \left[ |\overline{X}^\epsilon_1(s)| + |\overline{X}^\epsilon_2(s)|  \right]ds\\
&\quad + \Bigg| \int_0^t \left[ \left( C_{11} + D_1 \overline{F}_1(s) \right) \overline{X}^\epsilon_1(s) + \left( C_{12} + D_1 \overline{F}_2(s) \right) \overline{X}^\epsilon_2(s)  \right] dW(s)\Bigg|.
\end{align*}
Squaring and applying Ito's Isometry, we have that 
\begin{equation}\label{temp60}
\begin{split}
\E\left[ |\overline{X}^\epsilon_1(t)|^2 \right] &\leq 3 |x_1|^2 + K(T) \int_0^t \E \left[ |\overline{X}^\epsilon_1(s)|^2 + |\overline{X}^\epsilon_2(s)|^2  \right]ds\\
&\quad + \E \int_0^t \Big| \left( C_{11} + D_1 \overline{F}_1(s) \right) \overline{X}^\epsilon_1(s) + \left( C_{12} + D_1 \overline{F}_2(s) \right) \overline{X}^\epsilon_2(s)  \Big|^2 ds\\
&\leq 3 |x_1|^2 + K(T) \int_0^t \E \left[ |\overline{X}^\epsilon_1(s)|^2 + |\overline{X}^\epsilon_2(s)|^2  \right]ds.
\end{split}
\end{equation}
Now let us turn to the fast process \eqref{approxstate2} and fix $t \in [0,T]$. Applying Ito's formula to the mapping
\[[0,t] \ni s \mapsto e^{ (A_{22} + B_2 \overline{F}_2(t)) \frac{t-s}{\epsilon} } \overline{X}^\epsilon_2(s),
\] we obtain
\begin{align*}
d&\left( e^{ (A_{22} + B_2 \overline{F}_2(t)) \frac{t-s}{\epsilon} } \overline{X}^\epsilon_2(s) \right) \\
&= - \frac{1}{\epsilon}(A_{22} + B_2 \overline{F}_2(t)) e^{ (A_{22} + B_2 \overline{F}_2(t)) \frac{t-s}{\epsilon} } \overline{X}^\epsilon_2(s) ds + e^{ (A_{22} + B_2 \overline{F}_2(t)) \frac{t-s}{\epsilon} } d \overline{X}^\epsilon_2(s)\\
&=  \frac{1}{\epsilon} e^{ (A_{22} + B_2 \overline{F}_2(t)) \frac{t-s}{\epsilon} } \left[ \left( A_{21} + B_2 \overline{F}_1(s) \right) \overline{X}^\epsilon_1(s) + B_2 \left(\overline{F}_2(s) - \overline{F}_2(t) \right) \overline{X}_2^\epsilon(s) \right]ds\\
&\quad +  \frac{1}{\sqrt{\epsilon}}  e^{ (A_{22} + B_2 \overline{F}_2(t)) \frac{t-s}{\epsilon} }  \left[ \left( C_{21} + D_2 \overline{F}_1(s) \right) \overline{X}^\epsilon_1(s) + \left( C_{22} + D_2 \overline{F}_2(s) \right) \overline{X}^\epsilon_2(s)  \right] dW(s).
\end{align*} 
By another application of Ito's formula to the mapping
\[
[0,t] \ni s \mapsto |e^{ (A_{22} + B_2 \overline{F}_2(t)) \frac{t-s}{\epsilon} }  \overline{X}^\epsilon_2(s)|^2,
\] 
we obtain
\begin{align*}
d&|e^{ (A_{22} + B_2 \overline{F}_2(t)) \frac{t-s}{\epsilon} } \overline{X}^\epsilon_2(s)|^2\\
&= 2 \left\langle e^{ (A_{22} + B_2 \overline{F}_2(t)) \frac{t-s}{\epsilon} } \overline{X}^\epsilon_2(s), d\left( e^{ (A_{22} + B_2 \overline{F}_2(t)) \frac{t-s}{\epsilon} } \overline{X}^\epsilon_2(s)\right) \right\rangle\\
&\quad+ \left\langle d\left( e^{ (A_{22} + B_2 \overline{F}_2(t)) \frac{t-s}{\epsilon} } \overline{X}^\epsilon_2(s)\right), d\left( e^{ (A_{22} + B_2 \overline{F}_2(t)) \frac{t-s}{\epsilon} } \overline{X}^\epsilon_2(s)\right) \right\rangle\\
&= \frac{2}{\epsilon} \left\langle e^{ (A_{22} + B_2 \overline{F}_2(t)) \frac{t-s}{\epsilon} } \overline{X}^\epsilon_2(s),  e^{ (A_{22} + B_2 \overline{F}_2(t)) \frac{t-s}{\epsilon} }  \left[ \left( A_{21} + B_2 \overline{F}_1(s) \right) \overline{X}^\epsilon_1(s)  \right] \right\rangle ds\\
&\quad +\frac{2}{\epsilon} \left\langle e^{ (A_{22} + B_2 \overline{F}_2(t)) \frac{t-s}{\epsilon} } \overline{X}^\epsilon_2(s),  e^{ (A_{22} + B_2 \overline{F}_2(t)) \frac{t-s}{\epsilon} }  \left[ B_2 \left(\overline{F}_2(s) - \overline{F}_2(t) \right) \overline{X}^\epsilon_2(s) \right] \right\rangle ds\\
&\quad + \frac{2}{\sqrt{\epsilon}} \left\langle e^{ (A_{22} + B_2 \overline{F}_2(t)) \frac{t-s}{\epsilon} } \overline{X}^\epsilon_2(s) ,  e^{ (A_{22} + B_2 \overline{F}_2(t)) \frac{t-s}{\epsilon} }  \left( C_{21} + D_2 \overline{F}_1(s) \right) \overline{X}^\epsilon_1(s) dW(s)   \right\rangle \\
&\quad + \frac{2}{\sqrt{\epsilon}} \left\langle e^{ (A_{22} + B_2 \overline{F}_2(t)) \frac{t-s}{\epsilon} } \overline{X}^\epsilon_2(s) ,  e^{ (A_{22} + B_2 \overline{F}_2(t)) \frac{t-s}{\epsilon} }  \left( C_{22} + D_2 \overline{F}_2(s) \right) \overline{X}^\epsilon_2(s)   dW(s)   \right\rangle \\
&\quad + \frac{1}{\epsilon} \Big| e^{ (A_{22} + B_2 \overline{F}_2(t)) \frac{t-s}{\epsilon} }  \left[ \left( C_{21} + D_2 \overline{F}_1(s) \right) \overline{X}^\epsilon_1(s) + \left( C_{22} + D_2 \overline{F}_2(s) \right) \overline{X}^\epsilon_2(s)  \right] \Big|^2 ds.
\end{align*}
Integrating the above from $0$ to $t$ and taking the expectation gives
\begin{align*}
&\E \left[ |\overline{X}^\epsilon_2(t)|^2 \right] - |e^{ (A_{22} + B_2 \overline{F}_2(t) )\frac{t}{\epsilon} } x_2|^2\\
&= \frac{2}{\epsilon} \E \int_0^t \left\langle e^{ (A_{22} + B_2 \overline{F}_2(t)) \frac{t-s}{\epsilon} } \overline{X}^\epsilon_2(s),  e^{ (A_{22} + B_2 \overline{F}_2(t)) \frac{t-s}{\epsilon} }  \left[ \left( A_{21} + B_2 \overline{F}_1(s) \right) \overline{X}^\epsilon_1(s)  \right] \right\rangle ds\\
&\quad +\frac{2}{\epsilon} \E \int_0^t \left\langle e^{ (A_{22} + B_2 \overline{F}_2(t)) \frac{t-s}{\epsilon} } \overline{X}^\epsilon_2(s),  e^{ (A_{22} + B_2 \overline{F}_2(t)) \frac{t-s}{\epsilon} }  \left[ B_2 \left(\overline{F}_2(s) - \overline{F}_2(t) \right) \overline{X}^\epsilon_2(s) \right] \right\rangle ds\\
&\quad + \frac{2}{\sqrt{\epsilon}} \E \int_0^t \left\langle e^{ (A_{22} + B_2 \overline{F}_2(t)) \frac{t-s}{\epsilon} } \overline{X}^\epsilon_2(s) ,  e^{ (A_{22} + B_2 \overline{F}_2(t)) \frac{t-s}{\epsilon} }  \left( C_{21} + D_2 \overline{F}_1(s) \right) \overline{X}^\epsilon_1(s) dW(s)   \right\rangle \\
&\quad + \frac{2}{\sqrt{\epsilon}} \E \int_0^t \left\langle e^{ (A_{22} + B_2 \overline{F}_2(t)) \frac{t-s}{\epsilon} } \overline{X}^\epsilon_2(s) ,  e^{ (A_{22} + B_2 \overline{F}_2(t)) \frac{t-s}{\epsilon} }  \left( C_{22} + D_2 \overline{F}_2(s) \right) \overline{X}^\epsilon_2(s)   dW(s)   \right\rangle \\
&\quad + \frac{1}{\epsilon} \E \int_0^t \Big| e^{ (A_{22} + B_2 \overline{F}_2(t)) \frac{t-s}{\epsilon} }  \left[ \left( C_{21} + D_2 \overline{F}_1(s) \right) \overline{X}^\epsilon_1(s) + \left( C_{22} + D_2 \overline{F}_2(s) \right) \overline{X}^\epsilon_2(s)  \right] \Big|^2 ds.
\end{align*}
From Theorem 6.3 in Chapter 1 of \cite{yong1999stochastic},  for fixed $\epsilon \in (0,1]$, the processes $\overline{X}_1^\epsilon$ and $\overline{X}_2^\epsilon$ have bounded $4^{th}$ moments, and thus 
\begin{equation}\label{temp50}
\E \int_0^t \left\langle e^{ (A_{22} + B_2 \overline{F}_2(t)) \frac{t-s}{\epsilon} } \overline{X}^\epsilon_2(s) ,  e^{ (A_{22} + B_2 \overline{F}_2(t)) \frac{t-s}{\epsilon} }  \left( C_{21} + D_2 \overline{F}_1(s) \right) \overline{X}^\epsilon_1(s) dW(s)   \right\rangle = 0,
\end{equation}
and 
\begin{equation}\label{temp51}
\E\int_0^t \left\langle e^{ (A_{22} + B_2 \overline{F}_2(t)) \frac{t-s}{\epsilon} } \overline{X}^\epsilon_2(s) ,  e^{ (A_{22} + B_2 \overline{F}_2(t)) \frac{t-s}{\epsilon} }  \left( C_{22} + D_2 \overline{F}_2(s) \right) \overline{X}^\epsilon_2(s)   dW(s)   \right\rangle = 0.
\end{equation}
From Theorem \ref{theorem: DARE} and Corollary \ref{corollary: reducedsystem}, the eigenvalues of $A_{22} + B_2 \overline{F}_2(t)$ have negative real parts for all $t\in [0,T]$. Thus, there exists positive constants $M_\infty$ and $\gamma_\infty$ such that
\begin{equation}\label{temp52}
|e^{ (A_{22} + B_2 \overline{F}_2(t)) \frac{t-s}{\epsilon} }| \leq M_\infty e^{-\frac{\gamma_\infty (t-s)}{\epsilon}},\quad \forall 0 \leq s \leq t\leq T.
\end{equation}
Thus, using \eqref{temp50}-\eqref{temp52}, the Cauchy-Schwartz inequality and the uniform boundedness of $(\overline{F}_1(t),\overline{F}_2(t))$ for all $t\in [0,T]$, we obtain
\begin{align*}
\E \left[ |\overline{X}^\epsilon_2(t)|^2 \right] &\leq |e^{ (A_{22} + B_2 \overline{F}_2(t)) \frac{t}{\epsilon} } x_2|^2\\
&\quad+ \frac{2}{\epsilon}\E \int_0^t \left\langle e^{ (A_{22} + B_2 \overline{F}_2(t)) \frac{t-s}{\epsilon} } \overline{X}^\epsilon_2(s),  e^{ (A_{22} + B_2 \overline{F}_2(t)) \frac{t-s}{\epsilon} }  \left[ \left( A_{21} + B_2 \overline{F}_1(s) \right) \overline{X}^\epsilon_1(s)  \right] \right\rangle ds\\
&\quad +\frac{2}{\epsilon} \E \int_0^t \left\langle e^{ (A_{22} + B_2 \overline{F}_2(t)) \frac{t-s}{\epsilon} } \overline{X}^\epsilon_2(s),  e^{ (A_{22} + B_2 \overline{F}_2(t)) \frac{t-s}{\epsilon} }  \left[ B_2 \left(\overline{F}_2(s) - \overline{F}_2(t) \right) \overline{X}^\epsilon_2(s) \right] \right\rangle ds\\
&\quad + \frac{1}{\epsilon} \E \int_0^t \Big| e^{ (A_{22} + B_2 \overline{F}_2(t)) \frac{t-s}{\epsilon} }  \left[ \left( C_{21} + D_2 \overline{F}_1(s) \right) \overline{X}^\epsilon_1(s) + \left( C_{22} + D_2 \overline{F}_2(s) \right) \overline{X}^\epsilon_2(s)  \right] \Big|^2 ds\\
&\leq M_\infty^2 e^{-\frac{2\gamma_\infty t}{\epsilon}} |x_2|^2 + \frac{K(T)}{\epsilon} \int_0^t e^{-\frac{2\gamma_\infty (t-s)}{\epsilon}} \E \left[ |\overline{X}_1^\epsilon(s)|^2 + |\overline{X}_2^\epsilon(s)|^2 \right] ds\\
&\leq M_\infty^2 |x_2|^2 + \frac{K(T)}{\epsilon} \int_0^t e^{-\frac{2\gamma_\infty (t-s)}{\epsilon}} \E \left[ |\overline{X}_1^\epsilon(s)|^2 + |\overline{X}_2^\epsilon(s)|^2 \right] ds.
\end{align*}
Summing with \eqref{temp60}, we have
\begin{align*}
\E \left[ |\overline{X}^\epsilon_1(t)|^2 + |\overline{X}^\epsilon_2(t)|^2 \right] \leq 3 |x_1|^2 + M_\infty^2 |x_2|^2  + K(T) \int_0^t \left( 1 + \frac{1}{\epsilon}e^{-\frac{2\gamma_\infty (t-s)}{\epsilon}} \right) \E \left[ |\overline{X}^\epsilon_1(s)|^2 + |\overline{X}^\epsilon_2(s)|^2  \right]ds.
\end{align*}
By Gronwall's inequality (Theorem 15 of \cite{dragomir2003some}), we have that for fixed $\epsilon \in (0,1]$
\begin{align*}
\E \left[ |\overline{X}^\epsilon_1(t)|^2 + |\overline{X}^\epsilon_2(t)|^2 \right] &\leq \left(3 |x_1|^2 + M_\infty^2 |x_2|^2 \right) \exp\left[  K(T) \int_0^t \left(1 + \frac{1}{\epsilon} + \int_0^s \frac{-2\gamma_\infty}{\epsilon^2} e^{-\frac{2\gamma_\infty (s - r)}{\epsilon}} dr \right) ds \right]\\
&\leq \left(3 |x_1|^2 + M_\infty^2 |x_2|^2 \right) \exp\left[  K(T) \left( T + \frac{1}{2\gamma_\infty}  \right)\right].
\end{align*}
Since the right-hand side of the above inequality is independent of $\epsilon$ and $t$, we have that
\begin{align*}
\sup_{\epsilon \in (0,1]} \sup_{t\in [0,T]} \E \left[ |\overline{X}^\epsilon_1(t)|^2 + |\overline{X}^\epsilon_2(t)|^2 \right] \leq \left(3 |x_1|^2 + M_\infty^2 |x_2|^2 \right) \exp\left[  K(T) \left( T + \frac{1}{2\gamma_\infty}  \right)\right].
\end{align*}
\qed

The following theorem shows that using the approximately optimal control $\overline{u}^\epsilon$ defined in \eqref{approxoptimcontrol} gives a cost function close to the value function $V^\epsilon(x) = J^\epsilon(x; \widehat{u}^\epsilon)$ with an error of order $O(\epsilon)$.

\begin{theorem}
Suppose that Assumption \ref{mainassumptions} holds. Let $T > 0$ be any finite time horizon and $\epsilon^*$ be the small positive parameter defined in Theorem \ref{theorem: tikhonov}. Then, for every $x\in \R^n$, we have that for all $0 < \epsilon < \epsilon^*$
\begin{equation}
J^\epsilon(x; \overline{u}^\epsilon) - V^\epsilon(x) = O(\epsilon).
\end{equation}
\end{theorem}
\noindent \textit{Proof.} Let $P^\epsilon$ be the solution to the Riccati equation \eqref{riccati} and $X^\epsilon$ be the solution to the state equation \eqref{compactstates}. Applying Ito's formula to $\langle P^\epsilon(t) X^\epsilon(t), X^\epsilon(t) \rangle$ and by a completion of squares, we have that
\begin{equation}
J^\epsilon(x; u) = \frac{1}{2} \langle P^\epsilon(0) x,x\rangle - \frac{1}{2} \E \int_0^T \left[ | u(t) - \widehat{F}^\epsilon_1(t) X_1(t) - \widehat{F}^\epsilon_2(t) X^\epsilon_2(t)|^2 \right]dt
\end{equation}
From Theorem \ref{theorem: optimality}, $V^\epsilon(x) = \frac{1}{2} \langle P^\epsilon(0) x,x\rangle$. Hence we have that
\begin{align*}
|J^\epsilon(x; \overline{u}^\epsilon) - V^\epsilon(x)| &\leq \frac{1}{2} \E \int_0^T \left[ \Big| \left(\overline{F}_1(t) - \widehat{F}^\epsilon_1(t)\right) \overline{X}_1^\epsilon(t) + \left( \overline{F}_2(t) - \widehat{F}^\epsilon_2(t)\right) \overline{X}_2^\epsilon(t)\Big|^2 \right]dt\\
&\leq  \sup_{\epsilon \in (0,1]} \left\{ \left(\sup_{t\in [0,T]} \E\left[ |\overline{X}_1^\epsilon(t)|^2 \right]  \right) \int_0^T \left[ \Big| \overline{F}_1(t) - \widehat{F}^\epsilon_1(t) \Big|^2\right]dt\right\} \\
&\quad + \sup_{\epsilon \in (0,1]} \left\{ \left( \sup_{t\in [0,T]} \E\left[ |\overline{X}_2^\epsilon(t)|^2 \right]  \right) \int_0^T \left[ \Big| \overline{F}_2(t) - \widehat{F}^\epsilon_2(t) \Big|^2\right]dt \right\}.
\end{align*}
Hence Lemma \ref{lemma: feedbackbound} and Lemma \ref{lemma: statebound} gives the desired result.
\qed

Finally, we give an expression for the limiting value function.

\begin{theorem}\label{theorem: limitingvalue}
Suppose that Assumption \ref{mainassumptions} holds. Let $T > 0$ be any finite time horizon and $\epsilon^*$ be the small positive parameter defined in Theorem \ref{theorem: tikhonov}. Define
\begin{equation}
\overline{V}(x) = \frac{1}{2}\langle \overline{P}_{11}(0) x_1,x_1 \rangle.
\end{equation}
Then for all $0 < \epsilon < \epsilon^*$
\begin{equation}
V^\epsilon(x) - \overline{V}(x) = O(\epsilon). 
\end{equation}
\end{theorem}
\noindent \textit{Proof.} From Theorem \ref{theorem: optimality} and the first order partition \eqref{partition}, we have that
\begin{align*}
V^\epsilon(x) - \overline{V}(x) = \frac{1}{2} \langle P^\epsilon_{11}(0) x_1,x_1\rangle + \frac{\epsilon}{2} \left[ 2 \langle x_1, P^\epsilon_{12}(0)  x_2 \rangle + \langle P^\epsilon_{22}(0) x_2,x_2 \rangle \right] -  \frac{1}{2} \langle \overline{P}_{11}(0) x_1,x_1\rangle.
\end{align*}
From Theorem \ref{theorem: tikhonov}, we have that
\begin{align*}
|V^\epsilon(x) - \overline{V}(x)| &\leq \frac{1}{2} |\left\langle \left( P^\epsilon_{11}(0) - \overline{P}_{11}(0) \right) x_1,x_1\right\rangle| + \frac{\epsilon}{2} | 2 \langle P^\epsilon_{12}(0) x_1,x_2 \rangle + \langle P^\epsilon_{22}(0) x_2,x_2 \rangle | \\
&\leq \epsilon K(T, x)
\end{align*}
where $K(T,x)$ depends on $T$.
\qed

\bibliographystyle{abbrv}
\bibliography{references}

\end{document}